# TESTING CONVEX HYPOTHESES ON THE MEAN OF A GAUSSIAN VECTOR. APPLICATION TO TESTING QUALITATIVE HYPOTHESES ON A REGRESSION FUNCTION

By Yannick Baraud, Sylvie Huet and Béatrice Laurent

*Université de Nice Sophia-Antipolis, Laboratoire J. A. Dieudonné, INRA, Laboratoire de Biométrie and INSA de Toulouse, LSP*

In this paper we propose a general methodology, based on multiple testing, for testing that the mean of a Gaussian vector in $\mathbb{R}^n$ belongs to a convex set. We show that the test achieves its nominal level, and characterize a class of vectors over which the tests achieve a prescribed power. In the functional regression model this general methodology is applied to test some qualitative hypotheses on the regression function. For example, we test that the regression function is positive, increasing, convex, or more generally, satisfies a differential inequality. Uniform separation rates over classes of smooth functions are established and a comparison with other results in the literature is provided. A simulation study evaluates some of the procedures for testing monotonicity.

## 1. Introduction.

1.1. *The statistical framework.* We consider the following regression model:

$$Y_i = F(x_i) + \sigma\varepsilon_i, \qquad i = 1, \dots, n, \tag{1}$$

where $x_1 < x_2 < \cdots < x_n$ are known deterministic points in $[0,1]$, $\sigma$ is an unknown positive number and $(\varepsilon_i)_{i=1,\dots,n}$ is a sequence of i.i.d. unobservable standard Gaussian random variables. From the observation of $\mathbf{Y} = (Y_1, \dots, Y_n)'$, we consider the problem of testing that the regression function $F$ belongs to one of the following functional sets $\mathcal{K}$:

$$\mathcal{K}_{\geq 0} = \{F : [0,1] \to \mathbb{R},\ F \text{ is nonnegative}\}, \tag{2}$$

$$\mathcal{K}_{\nearrow} = \{F : [0,1] \to \mathbb{R},\ F \text{ is nondecreasing}\}, \tag{3}$$

$$\mathcal{K}_{\smile} = \{F : [0,1] \to \mathbb{R},\ F \text{ is nonconcave}\}, \tag{4}$$

$$\mathcal{K}_{r,R} = \left\{F : [0,1] \to \mathbb{R},\ \forall x \in [0,1],\ \frac{d^r}{dx^r}[R(x)F(x)] \geq 0\right\}. \tag{5}$$

In the above definition of $\mathcal{K}_{r,R}$, $r$ denotes a positive integer and $R$ a smooth, nonvanishing function from $[0,1]$ into $\mathbb{R}$. Choosing the function $R$ equal to 1 leads to test that the derivative of order $r$ is positive. Taking $r=1$ and choosing a suitable function $R$ leads to test that a positive function $F$ is decreasing at some prescribed rate. It is also possible to test that $F$ belongs to some classes of smooth functions. These testing hypotheses will be detailed in Section 3.

The problem is therefore to test some qualitative hypothesis on $F$. We shall show that it actually reduces to testing that the mean of the Gaussian vector $\mathbf{Y}$ belongs to a suitable convex subset of $\mathbb{R}^n$. Denoting by $\langle \cdot, \cdot \rangle$ the inner product of $\mathbb{R}^n$, this convex subset takes the form

$$\mathcal{C} = \{\mathbf{f} \in \mathbb{R}^n,\ \forall j \in \{1, \dots, p\}\ \langle \mathbf{f}, \mathbf{v}_j \rangle \leq 0\},$$

where the vectors $\{\mathbf{v}_1, \dots, \mathbf{v}_p\}$ are linearly independent in $\mathbb{R}^n$. The aim of this paper is to present a general methodology for the problem of testing that $\mathbf{f}$ belongs to $\mathcal{C}$ and to characterize a class of vectors over which the tests achieve a prescribed power. This general methodology is applied to test that the regression function $F$ belongs to one of the sets $\mathcal{K}$. For the procedures we propose, the least-favorable distribution under the null hypothesis is achieved for $F = 0$ and $\sigma = 1$. Consequently, by carrying out simulations, we easily obtain tests that achieve their nominal level for fixed values of $n$. Moreover, we show that these tests have good properties under smooth alternatives.

For the problem of testing positivity, monotonicity and convexity, we obtain tests based on the comparison of local means of consecutive observations. A precise description of these tests is given in Section 2. For the problem of testing monotonicity, our methodology also leads to tests based on the slopes of regression lines on short intervals, as explained in Section 3.1. These procedures, based on "running gradients," are akin to those proposed by Hall and Heckman (2000). For the problem of testing that $F$ belongs to $\mathcal{K}_{r,R}$ with a nonconstant function $R$ we refer the reader to Section 3.2. We have delayed the description of the general methodology for testing that $\mathbf{f}$ belongs to $\mathcal{C}$ to Section 4. Simulation studies for testing monotonicity are shown in Section 5. The proofs are postponed to Sections 6–9 and the Appendix.

1.2. *An overview of the literature.* In the literature tests of monotonicity have been widely studied in the regression model. The test proposed by Bowman, Jones and Gijbels (1998) is based on a procedure described in Silverman (1981) for testing unimodality of a density. This test is not powerful when the regression is flat or nearly flat, as emphasized by Hall and Heckman (2000). Hall and Heckman (2000) proposed a procedure based on "running gradients" over short intervals for which the least-favorable distribution under the null, when $\sigma$ is known, corresponds to the case where $F$ is identically constant. The test proposed by Gijbels, Hall, Jones and Koch (2000) is based on the signs of differences between observations. The test offers the advantage to not depend on the error distribution when it is continuous. Consequently, the nominal level of the test is guaranteed for all continuous error distributions. In the functional regression model with random $x_i$'s, the procedure proposed by Ghosal, Sen and van der Vaart (2000) is based on a locally weighted version of Kendall's tau. The procedure uses kernel smoothing with a particular choice of the bandwidth, and as in Gijbels, Hall, Jones and Koch (2000) depends on the signs of the quantities $(Y_j - Y_i)(x_i - x_j)$. They show that for certain local alternatives the power of their test tends to 1. Some comments on the power of our test under those alternatives can be found in Section 3.3. In Baraud, Huet and Laurent (2003b) we propose a procedure which aims at detecting discrepancies with respect to the $\mathbb{L}^2(\mu_n)$-distance where $\mu_n = n^{-1}\sum_{i=1}^n \delta_{x_i}$. This procedure generalizes that proposed in Baraud, Huet and Laurent (2003a) for linear hypotheses. A common feature of the present paper with these two lies in the fact that the proposed procedures achieve their nominal level and a prescribed power over a set of vectors we characterize. In the Gaussian white noise case, Juditsky and Nemirovski (2002) propose to test that the signal belongs to the cone of nonnegative, nondecreasing or nonconcave functions. For a given $r \in [1, +\infty[$, their tests are based on the estimation of the $\mathbb{L}^r$-distance between the signal and the cone. However, this approach requires that the signal have a known smoothness under the null. In the Gaussian white noise model, other tests of such qualitative hypotheses are proposed by Dümbgen and Spokoiny (2001). Their procedure is based on the supremum over all bandwidths of the distance in sup-norm between a kernel estimator and the null hypothesis. They adopt a minimax point of view to evaluate the performances of their tests and we adopt the same in Sections 2 and 3.

1.3. *Uniform separation rates and optimality.* Comparison of the performances of tests naturally arises in the problem of hypothesis testing. In this paper, we shall mainly describe the performances of our procedures in terms of uniform separation rates over classes of smooth functions. Given $\beta$ in $]0, 1[$, a class of smooth functions $\mathcal{F}$ and a "distance" $\Delta(\cdot)$ to the null hypothesis, we define the uniform separation rate of a test $\Phi$ over $\mathcal{F}$, denoted by $\rho(\Phi, \mathcal{F}, \Delta)$, as the smallest number $\rho$ such that the test guarantees a power not smaller than $1-\beta$ for all alternatives $F$ in $\mathcal{F}$ at distance $\rho$ from the null. More precisely,

$$\rho(\Phi,\mathcal{F},\Delta)=\inf\{\rho>0,\ \forall F\in\mathcal{F},\ \Delta(F)\geq\rho\Rightarrow\mathbb{P}_F(\Phi\text{ rejects })\geq 1-\beta\}. \tag{6}$$

In the regression or Gaussian white noise model, the word "rate" refers to the asymptotics of $\rho(\Phi,\mathcal{F},\Delta)=\rho_\tau(\Phi,\mathcal{F},\Delta)$ with respect to a scaling parameter $\tau$ (the number of observations $n$ in the regression model, the level of the noise in the Gaussian white noise). Comparing the performances of two tests of the same level amounts to comparing their uniform separation rates (the smaller the better). A test is said to be optimal if there exists no better test. The uniform separation rate of an optimal test is called the *minimax separation rate*. In the sequel, we shall enlarge this notion of optimality by saying that a test is rate-optimal over $\mathcal{F}$ if its uniform separation rate differs from the minimax one by a bounded function of $\tau$. Unfortunately, not much is known about the uniform separation rates of the tests mentioned in Section 1.2. The only exception we are aware of concerns the tests proposed by Dümbgen and Spokoiny (2001) and Juditsky and Nemirovski (2002) in the Gaussian white noise model (with $\tau=1/\sqrt{n}$), and Baraud, Huet and Laurent (2003b) in the regression model. The rates obtained by Juditsky and Nemirovski (2002) are established for the problem of testing that $F$ belongs to $\mathcal{K}\cap\mathcal{H}$, where $\mathcal{H}$ is a class of smooth functions. In contrast, in the papers by Baraud, Huet and Laurent (2003b) and Dümbgen and Spokoiny (2001), the null hypothesis is not restricted to those smooth functions belonging to $\mathcal{K}$. For the problem of testing positivity and monotonicity, Baraud, Huet and Laurent (2003b) established separation rates with respect to the $\mathbb{L}^2(\mu_n)$-distance to the null. For the problem of testing positivity, monotonicity and convexity, Dümbgen and Spokoiny (2001) considered the problem of detecting a discrepancy to the null in sup-norm. For any $L>0$, their procedures are proved to achieve the optimal rate $(L\log(n)/n)^{1/3}$ over the class of Lipschitz functions

$$\mathcal{H}_1(L)=\{F,\ \forall x,y\in[0,1],\ |F(x)-F(y)|\leq L|x-y|\}.$$

The optimality of this rate derives from the lower bounds established by Ingster [(1993), Section 2.4] for the more simple problem of testing $F=0$ against $F\neq 0$ in sup-norm. More generally, it can easily be derived from Ingster's results (see Proposition 2) that the minimax separation rate (in sup-norm) over Hölderian balls

$$\mathcal{H}_s(L)=\{F,\ \forall x,y\in[0,1],\ |F(x)-F(y)|\leq L|x-y|^s\}\quad\text{with } s\in\,]0,1] \tag{7}$$

is bounded from below (up to a constant) by $(L^{1/s}\log(n)/n)^{s/(1+2s)}$. In the regression setting, we propose tests of positivity, monotonicity and convexity whose uniform separation rates over $\mathcal{H}_s(L)$ achieve this lower bound whatever the value of $s \in ]\,0,1]$ and $L > 0$. In this paper, we discuss the optimality in the minimax sense over Hölderian balls with regularity $s$ in $]0,1]$ only. To our knowledge, the minimax rates over smoother classes of functions are unknown. It is beyond the scope of this paper to describe them.

For the problem of testing monotonicity or convexity, other choices of distance to the null are possible, for example, the distance in sup-norm between the first (resp. the second) derivative of $F$ and the set of nonnegative functions. For such choices, Dümbgen and Spokoiny also provided uniform separation rates for their tests. In the regression setting, the uniform separation rates we get coincide with their separation rates on the classes of functions they considered. We do not know whether these rates are optimal or not either in the Gaussian white noise model or in the regression model.

**2. Tests based on local means for testing positivity, monotonicity and convexity.** We consider the regression model given by (1) and propose tests of positivity, monotonicity and convexity for the function $F$. We first introduce some partitions of the design points and notation that will be used throughout the paper.

2.1. *Partition of the design points and notation.* We first define an almost regular partition of the set of indices $\{1,\dots,n\}$ into $\ell_n$ sets as follows: for each $k$ in $\{1,\dots,\ell_n\}$ we set

$$J_k = \left\{ i \in \{1,\dots,n\},\ \frac{k-1}{\ell_n} < \frac{i}{n} \leq \frac{k}{\ell_n} \right\}$$

and define the partition as

$$\mathcal{J}^{\ell_n} = \{J_k, k = \{1,\dots,\ell_n\}\}.$$

Then for each $\ell \in \{1,\dots,\ell_n\}$, we make a partition of $\{1,\dots,n\}$ into $\ell$ sets by gathering consecutive sets $J_k$. This partition is defined by

$$\mathcal{J}^{\ell} = \left\{ J_j^{\ell} = \bigcup_{(j-1)/\ell < k/\ell_n \leq j/\ell} J_k,\ j = 1,\dots,\ell \right\}. \tag{8}$$

We shall use the following notation.

(a) We use a bold type style for denoting the vectors of $\mathbb{R}^n$. We endow $\mathbb{R}^n$ with its Euclidean norm denoted by $\|\cdot\|$.

(b) For $\mathbf{v} \in \mathbb{R}^n$, let $\|\mathbf{v}\|_\infty = \max_{1\leq i\leq n} |v_i|$.

(c) For a linear subspace $V$ of $\mathbb{R}^n$, $\Pi_V$ denotes the orthogonal projector onto $V$.

(d) For $a \in \mathbb{R}_+$, $D \in \mathbb{N}\backslash\{0\}$ and $u \in [0,1]$, $\bar{\Phi}^{-1}(u)$ and $\bar{\chi}^{-1}_{D,a^2}(u)$ denote the $1-u$ quantile of, respectively, a standard Gaussian random variable and a noncentral $\chi^2$ with $D$ degrees of freedom and noncentrality parameter $a^2$.

(e) For $x \in \mathbb{R}$, $[x]$ denotes the integer part of $x$.

(f) For each $\mathbb{R}^n$-vector $\mathbf{v}$ and subset $J$ of $\{1,\dots,n\}$, we denote by $\mathbf{v}_J$ the $\mathbb{R}^n$-vector whose coordinates coincide with those of $\mathbf{v}$ on $J$ and vanish elsewhere. We denote by $\bar{\mathbf{v}}_J$ the quantity $\sum_{i\in J} v_i/|J|$.

(g) We denote by $\mathbf{1}$ the $\mathbb{R}^n$-vector $(1,\dots,1)'$ and by $\mathbf{e}_i$ the $i$th vector of the canonical basis.

(h) We define $V_{n,\mathrm{cste}}$ as the linear span of $\{\mathbf{1}_J, J \in \mathcal{J}^{\ell_n}\}$. Note that the dimension of $V_{n,\mathrm{cste}}$ equals $\ell_n$.

(i) The vector $\boldsymbol{\varepsilon}$ denotes a standard Gaussian variable in $\mathbb{R}^n$.

(j) We denote by $\mathbb{P}_{\mathbf{f},\sigma}$ the law of the Gaussian vector in $\mathbb{R}^n$ with expectation $\mathbf{f}$ and covariance matrix $\sigma^2 I_n$, where $I_n$ is the $n \times n$ identity matrix. We denote by $\mathbb{P}_{F,\sigma}$ the law of $Y$ under the model defined by (1).

(k) The level $\alpha$ of all our tests is chosen in $]0,1/2[$.

2.2. *Test of positivity.* We propose a level-$\alpha$ test for testing that $F$ belongs to $\mathcal{K}_{\geq 0}$ defined by (2). The testing procedure is based on the fact that if $F$ is nonnegative, then for any subset $J$ of $\{1,\dots,n\}$ the expectation of $\bar{\mathbf{Y}}_J$ is nonnegative. For $\ell \in \{1,\dots,\ell_n\}$, let $T_1^\ell(\mathbf{Y})$ be defined as

$$T_1^\ell(\mathbf{Y}) = \max_{J\in\mathcal{J}^\ell} \frac{-\sqrt{|J|}\,\bar{\mathbf{Y}}_J}{\|\mathbf{Y} - \Pi_{V_{n,\mathrm{cste}}}\mathbf{Y}\|}\sqrt{n-\ell_n},$$

and let $q_1(\ell,u)$ be the $1-u$ quantile of the random variable $T_1^\ell(\boldsymbol{\varepsilon})$. We introduce the test statistic

$$T_{\alpha,1} = \max_{\ell\in\{1,\dots,\ell_n\}} \{T_1^\ell(\mathbf{Y}) - q_1(\ell,u_\alpha)\}, \tag{9}$$

where $u_\alpha$ is defined as

$$u_\alpha = \sup\left\{u \in ]0,1[, \mathbb{P}\left(\max_{\ell\in\{1,\dots,\ell_n\}}\{T_1^\ell(\boldsymbol{\varepsilon}) - q_1(\ell,u)\} > 0\right) \leq \alpha\right\}. \tag{10}$$

We reject that $F$ belongs to $\mathcal{K}_{\geq 0}$ if $T_{\alpha,1}$ is positive.

COMMENT. When $\ell$ increases from 1 to $\ell_n$, the cardinality of the sets $J \in \mathcal{J}^\ell$ decreases. We thus take into account local discrepancies to the null hypothesis for various scales.

2.3. *Testing monotonicity.* We now consider the problem of testing that $F$ belongs to $\mathcal{K}_{\nearrow}$ defined by (3). The testing procedure relies on the following property: if $I$ and $J$ are two subsets of $\{1,2,\dots,n\}$ such that $I$ is on the left of $J$ and if $F \in \mathcal{K}_{\nearrow}$, then the expectation of the difference $\bar{\mathbf{Y}}_I - \bar{\mathbf{Y}}_J$ is nonpositive. For $\ell \in \{2, \ldots, \ell_n\}$, let $T_2^\ell(\mathbf{Y})$ be defined as

$$T_2^\ell(\mathbf{Y}) = \max_{1 \le i < j \le \ell} N_{ij}^\ell \frac{\bar{\mathbf{Y}}_{J_i^\ell} - \bar{\mathbf{Y}}_{J_j^\ell}}{\|\mathbf{Y} - \Pi_{V_{n,\mathrm{cste}}} \mathbf{Y}\|} \sqrt{n - \ell_n},$$

where

$$N_{ij}^\ell = \left( \frac{1}{|J_i^\ell|} + \frac{1}{|J_j^\ell|} \right)^{-1/2},$$

and let $q_2(\ell, u)$ be the $1 - u$ quantile of the random variable $T_2^\ell(\boldsymbol{\varepsilon})$. We introduce the test statistic

$$T_{\alpha,2} = \max_{\ell \in \{2,\ldots,\ell_n\}} \{T_2^\ell(\mathbf{Y}) - q_2(\ell, u)\}, \tag{11}$$

where $u_\alpha$ is defined as

$$u_\alpha = \sup\left\{ u \in ]0,1[, \mathbb{P}\left( \max_{\ell \in \{2,\ldots,\ell_n\}} \{T_2^\ell(\boldsymbol{\varepsilon}) - q_2(\ell, u)\} > 0 \right) \le \alpha \right\}. \tag{12}$$

We reject that $F$ belongs to $\mathcal{K}_{\nearrow}$ if $T_{\alpha,2}$ is positive.

2.4. *Testing convexity.* We now consider the problem of testing that $F$ belongs to $\mathcal{K}_{\smile}$ defined by (4). The testing procedure is based on the following property: if $I$, $J$ and $K$ are three subsets of $\{1, 2, \ldots, n\}$ such that $J$ is between $I$ and $K$ and if $F \in \mathcal{K}_{\smile}$, then we find a linear combination of $\bar{\mathbf{Y}}_I$, $\bar{\mathbf{Y}}_J$ and $\bar{\mathbf{Y}}_K$ with nonpositive expectation. Let $\mathbf{x} = (x_1, \ldots, x_n)'$ and for each $\ell \in \{3, \ldots, \ell_n\}$, $1 \le i < j < k \le \ell$, let

$$\lambda_{ijk}^\ell = \frac{\bar{\mathbf{x}}_{J_k^\ell} - \bar{\mathbf{x}}_{J_j^\ell}}{\bar{\mathbf{x}}_{J_k^\ell} - \bar{\mathbf{x}}_{J_i^\ell}}$$

and

$$N_{ijk}^\ell = \left( \frac{1}{|J_j^\ell|} + (\lambda_{ijk}^\ell)^2 \frac{1}{|J_i^\ell|} + (1 - \lambda_{ijk}^\ell)^2 \frac{1}{|J_k^\ell|} \right)^{-1/2}.$$

For $\ell \in \{3, \ldots, \ell_n\}$, let

$$T_3^\ell(\mathbf{Y}) = \max_{1 \le i < j < k \le \ell} N_{ijk}^\ell \frac{\bar{\mathbf{Y}}_{J_j^\ell} - \lambda_{ijk}^\ell \bar{\mathbf{Y}}_{J_i^\ell} - (1 - \lambda_{ijk}^\ell) \bar{\mathbf{Y}}_{J_k^\ell}}{\|\mathbf{Y} - \Pi_{V_{n,\mathrm{cste}}} \mathbf{Y}\| / \sqrt{n - \ell_n}},$$

and let $q_3(\ell, u)$ be the $1 - u$ quantile of the random variable $T_3^\ell(\boldsymbol{\varepsilon})$. We introduce the test statistic

$$T_{\alpha,3} = \max_{\ell \in \{3,\ldots,\ell_n\}} \{T_3^\ell(\mathbf{Y}) - q_3(\ell, u_\alpha)\}, \tag{13}$$

where $u_\alpha$ is defined as

$$u_\alpha = \sup\left\{u \in ]0,1[, \mathbb{P}\left(\max_{\ell\in\{3,\dots,\ell_n\}}\{T_3^\ell(\boldsymbol{\varepsilon}) - q_3(\ell,u)\} > 0\right) \leq \alpha\right\}. \tag{14}$$

We reject that $F$ belongs to $\mathcal{K}_\smile$ if $T_{\alpha,3}$ is positive.

2.5. *Properties of the procedures.* In this section we evaluate the performances of the previous procedures under the null and under smooth alternatives.

PROPOSITION 1. *Let $(T_\alpha, \mathcal{K})$ be either $(T_{\alpha,1}, \mathcal{K}_{\geq 0})$ or $(T_{\alpha,2}, \mathcal{K}_\nearrow)$ or $(T_{\alpha,3}, \mathcal{K}_\smile)$. We have*

$$\sup_{\sigma>0}\sup_{F\in\mathcal{K}} \mathbb{P}_{F,\sigma}(T_\alpha > 0) = \alpha.$$

*Assume now that $x_i = i/n$ for all $i = 1,\dots,n$ and $\ell_n = [n/2]$. Let us fix $\beta \in ]0,1[$ and define for each $s \in ]0,1]$ and $L > 0$*

$$\rho_n = L^{1/(1+2s)}\left(\frac{\sigma^2\log(n)}{n}\right)^{s/(1+2s)}.$$

*Then for $n$ large enough there exists some constant $\kappa$ depending on $\alpha, \beta, s$ only such that for all $F \in \mathcal{H}_s(L)$ satisfying*

$$\Delta(F) = \inf_{G\in\mathcal{K}} \|F - G\|_\infty \geq \kappa\rho_n \tag{15}$$

*we have*

$$\mathbb{P}_{F,\sigma}(T_\alpha > 0) \geq 1 - \beta.$$

COMMENT. This result states that our procedures are of size $\alpha$. Moreover, following the definition of the uniform separation rate of a test given in Section 1.3, this result shows that the tests achieve the uniform separation rate $\rho_n$ (in sup-norm) over the Hölderian ball $\mathcal{H}_s(L)$. In the following proposition, we show that this rate cannot be improved at least in the Gaussian white noise model for testing positivity and monotonicity. The proof can be extended to the case of testing convexity but is omitted here.

PROPOSITION 2. *Let $Y$ be the observation from the Gaussian white noise model*

$$dY(t) = F(t)\,dt + \frac{1}{\sqrt{n}}\,dW(t) \qquad \text{for } t \in [0,1], \tag{16}$$

*where $W$ is a standard Brownian motion. Let $\mathcal{K}$ be either the set $\mathcal{K}_{\geq 0}$ or $\mathcal{K}_\nearrow$ and let $\mathcal{F}$ be some class of functions. For the distance $\Delta(\cdot)$ to $\mathcal{K}$ given by (15), we define*

$$\rho_n(0,\mathcal{F}) = \inf \rho(\Phi, \mathcal{F}, \Delta),$$

*where $\rho(\Phi,\mathcal{F},\Delta)$ is given by (6) and where the infimum is taken over all tests $\Phi$ of level $3\alpha$ for testing "$F=0$." We define $\rho_n(\mathcal{K},\mathcal{F})$ similarly by taking the infimum over all tests $\Phi$ of level $\alpha$ for testing "$F\in\mathcal{K}$." The following inequalities hold*:

(i) *If $\mathcal{K}=\mathcal{K}_{\geq 0}$, then*

$$\rho_n(\mathcal{K},\mathcal{F})\geq\rho_n(0,\mathcal{F}).$$

*If $\mathcal{K}=\mathcal{K}_{\nearrow}$, then for some constant $\kappa$ depending on $\alpha$ and $\beta$ only*

$$\rho_n(\mathcal{K},\mathcal{F})\geq\frac{1}{2}\left[\rho_n(0,\mathcal{F})-\kappa\frac{\sigma}{\sqrt{n}}\right].$$

(ii) *In particular, if $\mathcal{F}=\mathcal{H}_s(L)$, for $n$ large enough there exists some constant $\kappa'$ depending on $\alpha,\beta$ and $s$ only such that*

$$\rho_n(\mathcal{K},\mathcal{F})\geq\kappa' L^{1/(1+2s)}\left(\frac{\log(n)}{n}\right)^{s/(1+2s)}. \tag{17}$$

The proof of the first part of the proposition extends easily to the regression framework. The second part (ii), namely (17), derives from (i) and the lower bound on $\rho_n(0,\mathcal{F})$ established by Ingster (1993).

For the problem of testing the positivity of a signal in the Gaussian white noise model, Juditsky and Nemirovski (2002) showed that the minimax separation rate with respect to the $\mathbb{L}^r$-distance ($r\in[1,+\infty[$) is of the same order as $\rho_n$ up to a logarithmic factor.

## 3. Testing that $F$ satisfies a differential inequality.

In this section, we consider the problem of testing that $F$ belongs to $\mathcal{K}_{r,R}$ defined by (5). Several applications of such hypotheses can be of interest. For example, by taking $r=1$ and $R(x)=-\exp(ax)$ (for some positive number $a$), one can test that a positive function $F$ is decreasing at rate $\exp(-ax)$, that is, satisfies

$$\forall\, x\in[0,1]\qquad 0<F(x)\leq F(0)\exp(-ax).$$

Other kinds of decay are possible by suitably choosing the function $R$. Another application is to test that $F$ belongs to the class of smooth functions

$$\{F:[0,1]\to\mathbb{R},\ \|F^{(r)}\|_\infty\leq L\}.$$

To tackle this problem, it is enough to test that the derivatives of order $r$ of the functions $F_1(x)=-F(x)+Lx^r/r!$ and $F_2(x)=F(x)+Lx^r/r!$ are positive. This is easily done by considering a multiple testing procedure based on the data $-Y_i+Lx_i^r/r!$ for testing that $F_1$ is positive, and on $Y_i+Lx_i^r/r!$ for testing that $F_2$ is positive.

In Section 3.1 we consider the case where the function $R$ equals 1. The procedure then amounts to testing that the derivative of order $r$ of $F$ is nonnegative. We turn to the general case in Section 3.2.

We first introduce the following notation.

(a) For $\mathbf{w} \in \mathbb{R}^n$, we denote by $R \star \mathbf{w}$ the vector whose $i$th coordinate $(R \star \mathbf{w})_i$ equals $R(x_i)w_i$.

(b) For $k \in \mathbb{N} \setminus \{0\}$, we denote by $\mathbf{w}^k$ the $\mathbb{R}^n$-vector $(w_1^k, \dots, w_n^k)$, and we set $\mathbf{w}^0 = \mathbf{1}$ by convention.

(c) For $J \subset \{1, \dots, n\}$, let us define $\mathcal{X}_J$ as the space spanned by $\mathbf{1}_J, \mathbf{x}_J, \dots, \mathbf{x}_J^{r-1}$.

3.1. *Testing that the derivative of order $r$ of $F$ is nonnegative.* In this section we take $R(x) = 1$ for all $x \in [0,1]$. The procedure relies on the idea that if the derivative of order $r$ of $F$ is nonnegative, then on each subset $J$ of $\{1, 2, \dots, n\}$, the highest degree coefficient of the polynomial regression of degree $r$ based on the pairs $\{(x_i, F(x_i)), i \in J\}$ is nonnegative. For example, under the assumption that $F$ is nondecreasing, the slope of the regression based on the pairs $\{(x_i, F(x_i)), i \in J\}$ is nonnegative.

Let $\ell_n = [n/(2(r+1))]$, let $V_n$ be the linear span of $\{\mathbf{1}_J, \mathbf{x}_J, \dots, \mathbf{x}_J^r,\ J \in \mathcal{J}^{\ell_n}\}$, and for each $J \subset \{1, \dots, n\}$

$$\mathbf{t}_J^* = -\frac{\mathbf{x}_J^r - \Pi_{\mathcal{X}_J}\mathbf{x}_J^r}{\|\mathbf{x}_J^r - \Pi_{\mathcal{X}_J}\mathbf{x}_J^r\|}.$$

For each $\ell \in \{1, \dots, \ell_n\}$, let $T^\ell(\mathbf{Y})$ be defined as

$$T^\ell(\mathbf{Y}) = \max_{J \in \mathcal{J}^\ell} \frac{\langle \mathbf{Y}, \mathbf{t}_J^* \rangle}{\|\mathbf{Y} - \Pi_{V_n}\mathbf{Y}\|}\sqrt{n - d_n} \tag{18}$$

and let $q(\ell, u)$ denote the $1-u$ quantile of the random variable $T^\ell(\boldsymbol{\varepsilon})$. We introduce the following test statistic:

$$T_\alpha = \max_{\ell \in \{1, \dots, \ell_n\}} \{T^\ell(\mathbf{Y}) - q(\ell, u_\alpha)\}, \tag{19}$$

where $u_\alpha$ is defined as

$$u_\alpha = \sup\left\{u \in ]\,0,1[, \mathbb{P}\left(\max_{\ell \in \{1, \dots, \ell_n\}} \{T^\ell(\boldsymbol{\varepsilon}) - q(\ell, u)\} > 0\right) \le \alpha\right\}. \tag{20}$$

We reject the null hypothesis if $T_\alpha$ is positive.

COMMENT. When $r = 1$, the procedure is akin to that proposed by Hall and Heckman (2000) where for all $\ell$, $q(\ell, u_\alpha)$ is the $1-\alpha$ quantile of $\max_{\ell \in \{1, \dots, \ell_n\}} T^\ell(\boldsymbol{\varepsilon})$.

3.2. *Extension to the general case.* The ideas underlying the preceding procedures extend to the case where $R \not\equiv 1$. In the general case, the test is obtained as follows.

Let $\ell_n$ be such that the dimension $d_n$ of the linear space

$$V_n = \mathrm{Span}\{\mathbf{1}_J, \mathbf{x}_J, \dots, \mathbf{x}_J^r, R \star \mathbf{1}_J, \dots, R \star \mathbf{x}_J^r, J \in \mathcal{J}^{\ell_n}\} \tag{21}$$

is not larger than $n/2$. We define for each $J \subset \{1, \ldots, n\}$

$$\mathbf{t}_J^* = -\frac{R \star (\mathbf{x}_J^r - \Pi_{\mathcal{X}_J}\mathbf{x}_J^r)}{\gamma_J} \qquad \text{where } \gamma_J = \|R \star (\mathbf{x}_J^r - \Pi_{\mathcal{X}_J}\mathbf{x}_J^r)\|. \tag{22}$$

We reject that $F$ belongs to $\mathcal{K}_{r,R}$ if $T_\alpha$ defined by (19) is positive.

3.3. *Properties of the tests.* In this section we describe the behavior of the procedure. We start with some notation.

(a) Let us define the function $\Lambda(F)$ as

$$\Lambda(F)(x) = \frac{d^r}{dx^r}[R(x)F(x)],$$

and let $\omega$ be its modulus of continuity defined for all $h > 0$ by

$$\omega(h) = \sup_{|x-y| \le h} |\Lambda(F)(x) - \Lambda(F)(y)|.$$

(b) For $J \in \bigcup_{\ell=1}^{\ell_n} \mathcal{J}^\ell$, let us denote by $x_J^-$ (resp. $x_J^+$) the quantities $\min\{x_i,\ i \in J\}$ (resp. $\max\{x_i,\ i \in J\}$) and set $h_J = x_J^+ - x_J^-$.

(c) Let $\mathbf{f} = (F(x_1), \ldots, F(x_n))'$ and for each $\ell = 1, \ldots, \ell_n$ and $\beta \in ]0,1[$, let

$$\nu_\ell(\mathbf{f}, \beta) = \left(\frac{q(\ell, u_\alpha)}{\sqrt{n - d_n}}\sqrt{\bar{\chi}^{-1}_{n-d_n, \|\mathbf{f} - \Pi_{V_n}\mathbf{f}\|^2/\sigma^2}(\beta/2)} + \bar{\Phi}^{-1}(\beta/2)\right)\sigma. \tag{23}$$

(d) For each $\rho > 0$, let

$$\mathcal{E}_{n,r}(\rho) = \left\{F : [0,1] \to \mathbb{R},\ F^{(r)} \in \mathcal{H}_s(L), -\inf_{x \in [0,1]} F^{(r)}(x) \ge \rho\right\}.$$

We have the following result.

Proposition 3. *Let $T_\alpha$ be the test statistic defined in Section 3.2. We have*

$$\sup_{\sigma > 0} \sup_{F \in \mathcal{K}_{r,R}} \mathbb{P}_{F,\sigma}(T_\alpha > 0) = \alpha.$$

*For each $\beta \in ]0,1[$ we have*

$$\mathbb{P}_{F,\sigma}(T_\alpha > 0) \ge 1 - \beta,$$

*if for some $\ell \in \{1, \ldots, \ell_n\}$ there exists a set $J \in \mathcal{J}^\ell$ such that either*

$$-\inf_{i \in J} \Lambda(F)(x_i) \ge \nu_\ell(\mathbf{f}, \beta)\frac{r!\gamma_J}{\|\mathbf{x}_J^r - \Pi_{\mathcal{X}_J}\mathbf{x}_J^r\|^2} + \omega(h_J), \tag{24}$$

*or*

$$\inf_{x \in ]x_J^-, x_J^+[} -\Lambda(F)(x) \ge \nu_\ell(\mathbf{f}, \beta)\frac{r!\gamma_J}{\|\mathbf{x}_J^r - \Pi_{\mathcal{X}_J}\mathbf{x}_J^r\|^2}. \tag{25}$$

*Moreover, if $R \equiv 1$, then there exists some constant $\kappa$ depending on $\alpha, \beta, s$ and $r$ only such that for $n$ large enough and for all $F \in \mathcal{E}_{n,r}(\rho_{n,r})$ with*

$$\rho_{n,r} = \kappa\left(\frac{\sigma^2 \log(n)}{n}\right)^{s/(1+2(s+r))} L^{(1+2r)/(1+2(s+r))}$$

*we have*

$$\mathbb{P}_{F,\sigma}(T_\alpha > 0) \geq 1 - \beta.$$

COMMENT 1. In the particular case where $R \equiv 1$, let us give the orders of magnitude of the quantities appearing in the above proposition. Under the assumption that $\|\mathbf{f} - \Pi_{V_n}\mathbf{f}\|^2/n$ is smaller than $\sigma^2$, one can show that $\nu_\ell$ is of order $\sqrt{\log(n)}$ (see Section 9.2). When $R \equiv 1$, we have $\gamma_J = \|\mathbf{x}_J^r - \Pi_{\mathcal{X}_J}\mathbf{x}_J^r\|$ and it follows from computations that will be detailed in the proofs that

$$\nu_\ell(\mathbf{f}, \beta)\frac{r!\gamma_J}{\|\mathbf{x}_J^r - \Pi_{\mathcal{X}_J}\mathbf{x}_J^r\|^2} \leq C\sqrt{\frac{\log(n)}{nh_J^{1+2r}}} \tag{26}$$

for some constant $C$ which does not depend on $J$ or $n$.

COMMENT 2. In the particular case where $r = 1$, (26) allows us to compare our result to the performance of the test proposed by Ghosal, Sen and van der Vaart (2000). For each $\delta \in ]0, 1/3[$, they give a procedure (depending on $\delta$) that is powerful if the function $F$ is continuously differentiable and satisfies that for all $x$ in some interval of length $n^{-\delta}$, $F'(x) < -M\sqrt{\log(n)}n^{-(1-3\delta)/2}$ for some $M$ large enough.

By using (25) and the upper bound in (26) with $h_J$ of order $n^{-\delta}$, we deduce from Proposition 3 that our procedure is powerful too over this class of functions. Note that by considering a multiple testing procedure based on various scales $\ell$, our test does not depend on $\delta$ and is therefore powerful for all $\delta$ simultaneously.

COMMENT 3. For $r = 1$ (resp. $r = 2$) and $s = 1$, Dümbgen and Spokoiny (2001) obtained the uniform separation rate $\rho_{n,r}$ for testing monotonicity (resp. convexity) in the Gaussian white noise model.

COMMENT 4. For the problem of testing monotonicity ($r = 1$ and $R \equiv 1$), it is possible to combine this procedure with that proposed in Section 2.3. More precisely, consider the test which rejects the null at level $2\alpha$ if one of these two tests rejects. The so-defined test performs as well as the best of these two tests under the alternative.

**4. A general approach.** The problems we have considered previously reduce to testing that $\mathbf{f} = (F(x_1), \dots, F(x_n))'$ belongs to a convex set of the form

$$\mathcal{C} = \{\mathbf{f} \in \mathbb{R}^n,\ \forall j \in \{1, \dots, p\}\ \langle \mathbf{f}, \mathbf{v}_j \rangle \le 0\}, \tag{27}$$

where the vectors $\{\mathbf{v}_1, \dots, \mathbf{v}_p\}$ are linearly independent in $\mathbb{R}^n$. For example, testing that the regression function $F$ is nonnegative or nondecreasing amounts to testing that the mean of $\mathbf{Y}$ belongs, respectively, to the convex subsets of $\mathbb{R}^n$

$$\mathcal{C}_{\ge 0} = \{\mathbf{f} \in \mathbb{R}^n,\ \forall i \in \{1, \dots, n\}\ f_i \ge 0\} \tag{28}$$

and

$$\mathcal{C}_{\nearrow} = \{\mathbf{f} \in \mathbb{R}^n,\ \forall i \in \{1, \dots, n-1\}\ f_{i+1} - f_i \ge 0\}. \tag{29}$$

Clearly, these sets are of the form given by (27) by taking, respectively, $p = n$, $\mathbf{v}_j = -\mathbf{e}_j$ and $p = n-1$, $\mathbf{v}_j = \mathbf{e}_j - \mathbf{e}_{j+1}$. The following proposition extends this result to the general case. Note that one can also define the set $\mathcal{C}$ as

$$\mathcal{C} = \{\mathbf{f} \in \mathbb{R}^n,\ L_1(\mathbf{f}) \ge 0, \dots, L_p(\mathbf{f}) \ge 0\},$$

where the $L_i$'s are $p$ independent linear forms. We shall use this definition of $\mathcal{C}$ in the following.

PROPOSITION 4. *For each $r \in \{1, \dots, n-1\}$ and $i \in \{1, \dots, n-r\}$ let $\phi_{i,r}$ be the linear form defined for $\mathbf{w} \in \mathbb{R}^n$ by*

$$\phi_{i,r}(\mathbf{w}) = \det \begin{pmatrix} 1 & x_i & \cdots & x_i^{r-1} & w_i \\ 1 & x_{i+1} & \cdots & x_{i+1}^{r-1} & w_{i+1} \\ \vdots & \vdots & \vdots & \vdots & \vdots \\ 1 & x_{i+r} & \cdots & x_{i+r}^{r-1} & w_{i+r} \end{pmatrix}.$$

*If $F$ belongs to $\mathcal{K}_{\smile}$, then $\mathbf{f} = (F(x_1), \dots, F(x_n))'$ belongs to*

$$\mathcal{C}_{\smile} = \{\mathbf{f} \in \mathbb{R}^n,\ \forall i \in \{1, \dots, n-2\},\ \phi_{i,2}(\mathbf{f}) \ge 0\}. \tag{30}$$

*If $F$ belongs to $\mathcal{K}_{r,R}$, then $\mathbf{f}$ belongs to*

$$\mathcal{C}_{r,R} = \{\mathbf{f} \in \mathbb{R}^n,\ \forall i \in \{1, \dots, n-r\},\ \phi_{i,r}(R \star \mathbf{f}) \ge 0\}.$$

With the aim of keeping our notation as simple as possible, we omit the dependence of the linear forms $\phi_{i,r}$ on $r$ when there is no ambiguity. The remaining part of the section is organized as follows. In the next section we present a general approach for the problem of testing that $\mathbf{f}$ belongs to $\mathcal{C}$. In the last section we show how this approach applies to the problems of hypothesis testing considered in Sections 2 and 3.

4.1. *Testing that* $\mathbf{f}$ *belongs to* $\mathcal{C}$. We consider the problem of testing that the vector $\mathbf{f} = (f_1, \ldots, f_n)'$ involved in the regression model

$$Y_i = f_i + \sigma\varepsilon_i, \qquad i = 1, \ldots, n, \tag{31}$$

belongs to $\mathcal{C}$ defined by (27). Our aim is twofold: first, build a test which achieves its nominal level, and second, describe for each $n$ a class of vectors over which this test is powerful.

*The testing procedure.* The testing procedure relies on the following idea: since under the assumption that $\mathbf{f}$ belongs to $\mathcal{C}$, the quantities $\langle \mathbf{f}, \sum_{j=1}^p \lambda_j \mathbf{v}_j \rangle$ are nonpositive for all nonnegative numbers $\lambda_1, \ldots, \lambda_p$ we base our test statistic on random variables of the form $\langle \mathbf{Y}, \sum_{j=1}^p \lambda_j \mathbf{v}_j \rangle$ for nonnegative sequences of $\lambda_j$'s.

We denote by $\mathcal{T}$ the subset of $\mathbb{R}^n$ defined by

$$\mathcal{T} = \left\{ \mathbf{t} = \sum_{j=1}^p \lambda_j \mathbf{v}_j,\ \|\mathbf{t}\| = 1,\ \lambda_j \geq 0,\ \forall j = 1, \ldots, p \right\}. \tag{32}$$

Let $\mathcal{T}_n$ be a finite subset of $\mathcal{T}$ such that there exists some linear space $V_n$ with dimension $d_n < n$ containing the linear span of $\mathcal{T}_n$. Let $\{q_{\mathbf{t}}(\alpha), \mathbf{t} \in \mathcal{T}_n\}$ be a sequence of numbers satisfying

$$\mathbb{P}\left[ \sup_{\mathbf{t} \in \mathcal{T}_n} \left( \sqrt{n - d_n} \frac{\langle \boldsymbol{\varepsilon}, \mathbf{t} \rangle}{\|\boldsymbol{\varepsilon} - \Pi_{V_n} \boldsymbol{\varepsilon}\|} - q_{\mathbf{t}}(\alpha) \right) > 0 \right] = \alpha. \tag{33}$$

We reject the null hypothesis if the statistic

$$T_\alpha = \sup_{\mathbf{t} \in \mathcal{T}_n} \left( \sqrt{n - d_n} \frac{\langle \mathbf{Y}, \mathbf{t} \rangle}{\|\mathbf{Y} - \Pi_{V_n} \mathbf{Y}\|} - q_{\mathbf{t}}(\alpha) \right) \tag{34}$$

is positive.

*Properties of the test.* For all $\beta \in ]0, 1[$ and each $\mathbf{t} \in \mathcal{T}_n$ let

$$v_{\mathbf{t}}(\mathbf{f}, \beta) = \left( q_{\mathbf{t}}(\alpha) \frac{1}{\sqrt{n - d_n}} \sqrt{\bar{\chi}^{-1}_{n-d_n, \|\mathbf{f} - \Pi_{V_n} \mathbf{f}\|^2/\sigma^2}(\beta/2)} + \bar{\Phi}^{-1}(\beta/2) \right) \sigma. \tag{35}$$

The order of magnitude of $v_{\mathbf{t}}(\mathbf{f}, \beta)$ is proved to be $\sqrt{\log(n)}\,\sigma$ under the assumption that $\|f - \Pi_{V_n} f\|^2/n$ is smaller than $\sigma^2$ as is shown in the proof of Proposition 1.

We have the following result.

THEOREM 1. *Let $T_\alpha$ be the test statistic defined by (34). We have*

$$\sup_{\sigma > 0} \sup_{\mathbf{f} \in \mathcal{C}} \mathbb{P}_{\mathbf{f}, \sigma}(T_\alpha > 0) = \mathbb{P}_{0,1}(T_\alpha > 0) = \alpha. \tag{36}$$

*Moreover, if there exists $\mathbf{t} \in \mathcal{T}_n$ such that $\langle \mathbf{f}, \mathbf{t} \rangle \geq v_{\mathbf{t}}(\mathbf{f}, \beta)$, then*

$$\mathbb{P}_{\mathbf{f}, \sigma}(T_\alpha > 0) \geq 1 - \beta.$$

COMMENT. The values of the $q_{\mathbf{t}}(\alpha)$'s that satisfy (33) can be easily obtained by simulations under $\mathbb{P}_{0,1}$. This property of our procedure lies in the fact that the least-favorable distribution under the null is $\mathbb{P}_{0,1}$. Note that we do not need to use bootstrap procedures to implement the test.

4.2. *How to apply these procedures to test qualitative hypotheses.* In the sequel, we give the choices of $\mathcal{T}_n$ and $V_n$ leading to the tests presented in Sections 2 and 3.

*For the test of positivity described in Section 2.2.* We take $\mathcal{T}_n = \mathcal{T}_{n,1}$, with $\mathcal{T}_{n,1} = \bigcup_{\ell=1}^{\ell_n} \mathcal{T}_{n,1}^\ell$, where for all $\ell \in \{1,\dots,\ell_n\}$

$$\mathcal{T}_{n,1}^\ell = \left\{ -\frac{1}{\sqrt{|J|}} \sum_{j\in J} \mathbf{e}_j,\ J \in \mathcal{J}^\ell \right\}.$$

We take $V_n = V_{n,\mathrm{cste}}$. Note that $V_{n,\mathrm{cste}}$ is also the linear span of $\mathcal{T}_{n,1}$.

*For the test of monotonicity described in Section 2.3.* Let us define for each $\ell \in \{2,\dots,\ell_n\}$ and $1 \le i < j \le \ell$,

$$\mathbf{e}_{ij}^\ell = N_{ij}^\ell \left( \frac{1}{|J_i^\ell|} \sum_{l\in J_i^\ell} \mathbf{e}_l - \frac{1}{|J_j^\ell|} \sum_{l \in J_j^\ell} \mathbf{e}_l \right). \tag{37}$$

Note that $N_{ij}^\ell$ is such that $\|\mathbf{e}_{ij}^\ell\| = 1$. We take $\mathcal{T}_n = \mathcal{T}_{n,2}$, with $\mathcal{T}_{n,2} = \bigcup_{\ell=2}^{\ell_n} \mathcal{T}_{n,2}^\ell$, where

$$\mathcal{T}_{n,2}^\ell = \{\mathbf{e}_{ij}^\ell,\ 1 \le i < j \le \ell\},$$

and we take $V_n = V_{n,\mathrm{cste}}$. Note that $V_n$ contains $\mathcal{T}_{n,2}$.

*For the test of convexity presented in Section 2.4.* Let us define for each $\ell \in \{3,\dots,\ell_n\}$, $1 \le i < j < k \le \ell$,

$$\mathbf{e}_{ijk}^\ell = N_{ijk}^\ell \left( \frac{1}{|J_j^\ell|} \sum_{l\in J_j^\ell} \mathbf{e}_l - \lambda_{ijk}^\ell \frac{1}{|J_i^\ell|} \sum_{l\in J_i^\ell} \mathbf{e}_l - (1-\lambda_{ijk}^\ell) \frac{1}{|J_k^\ell|} \sum_{l\in J_k^\ell} \mathbf{e}_l \right). \tag{38}$$

Note that $N_{ijk}^\ell$ is such that $\|\mathbf{e}_{ijk}^\ell\| = 1$. We take $\mathcal{T}_n = \mathcal{T}_{n,3}$, with $\mathcal{T}_{n,3} = \bigcup_{\ell=3}^{\ell_n} \mathcal{T}_{n,3}^\ell$, where

$$\mathcal{T}_{n,3}^\ell = \{\mathbf{e}_{ijk}^\ell,\ 1 \le i < j < k \le \ell\},$$

and we take $V_n = V_{n,\mathrm{cste}}$. Note that $V_n$ contains $\mathcal{T}_{n,3}$.

*For the test of* $F \in \mathcal{K}_{r,R}$ *presented in Section* 3. We take

$$\mathcal{T}_{n,4} = \bigcup_{\ell=1}^{\ell_n} \mathcal{T}_n^\ell \qquad \text{where } \mathcal{T}_n^\ell = \{\mathbf{t}_J^*,\ J \in \mathcal{J}^\ell\}$$

and $V_n = V_{n,4}$ defined by (21). Note that $V_n$ contains $\mathcal{T}_{n,4}$.

We justify these choices of $\mathcal{T}_n$ by the following proposition proved in Section 7.

PROPOSITION 5. *Let* $\mathcal{C}$ *and* $\mathcal{T}_n$ *be either* $(\mathcal{C}_{\geq 0}, \mathcal{T}_{n,1})$, $(\mathcal{C}_{\nearrow}, \mathcal{T}_{n,2})$, $(\mathcal{C}_{\smile}, \mathcal{T}_{n,3})$ *or* $(\mathcal{C}_{r,R}, \mathcal{T}_{n,4})$. *There exist* $\mathbf{v}_1, \dots, \mathbf{v}_p$ *for which* $\mathcal{C}$ *is of the form given by* (27) *and for which* $\mathcal{T}$ *defined by* (32) *contains* $\mathcal{T}_n$.

**5. Simulation studies.** In this section we describe how to implement the test for $F \in \mathcal{K}_{\nearrow}$ and we carry out a simulation study in order to evaluate the performances of our tests both when the errors are Gaussian and when they are not. We first describe how the testing procedure is performed, then we present the simulation experiment and finally discuss the results of the simulation study.

5.1. *The testing procedures.* We carry out the simulation study for the two testing procedures described in Sections 2.3 and 3.1. In the sequel, the procedure based on differences of *local means* and described in Section 2.3 is called LM and the procedure based on *local gradients* defined below (from the test statistic given in Section 3.1 with $r = 1$) is called LG.

In the case of the procedure LM, we set $T_{\mathrm{LM}} = T_{\alpha,2}$ defined in (11). For each $\ell$, the quantiles $q_2(\ell, u_\alpha)$ are calculated as follows. For $u$ varying among a suitable grid of values $u_1, \dots, u_m$, we estimate by simulations the quantity

$$p(u_j) = \mathbb{P}\left(\max_{l=1,\dots,\ell_n} \{T_2^\ell(\boldsymbol{\varepsilon}) - q_2(\ell, u_j)\} > 0\right),$$

$\boldsymbol{\varepsilon}$ being an $n$-sample of $\mathcal{N}(0,1)$, and we take $u_\alpha$ as $\max\{u_j, p(u_j) \leq \alpha\}$. Note that $u_\alpha$ does not depend on $(x_i, i = 1, \dots, n)$, but only on the number of observations $n$.

In the case of the procedure LG, the test statistic is defined as follows. For each $\ell = 1, \dots, \ell_n$ and for $J \in \mathcal{J}^\ell$, we take

$$\mathbf{t}_J^* = \frac{\bar{x}_J \mathbf{1}_J - \mathbf{x}_J}{\|\bar{x}_J \mathbf{1}_J - \mathbf{x}_J\|}.$$

The space $V_n$ reduces to $V_{n,\mathrm{lin}}$, the linear space of dimension $2\ell_n$ generated by

$$\{\mathbf{1}_J, \mathbf{x}_J, J \in \mathcal{J}^{\ell_n}\}.$$

The test statistic $T_\alpha$ takes the form

$$T_{\mathrm{LG}} = T_{\alpha,4} = \max_{\ell=1,\dots,\ell_n} \{T_4^\ell(\mathbf{Y}) - q_4(\ell, u_\alpha)\},$$

where for each $\ell \in \{1, \dots, \ell_n\}$,

$$T_4^\ell(\mathbf{Y}) = \max_{J\in\mathcal{J}^\ell} \sqrt{n-2\ell_n}\frac{\langle \mathbf{Y}, \mathbf{t}_J^*\rangle}{\|\mathbf{Y} - \Pi_{V_{n,\mathrm{lin}}}\mathbf{Y}\|},$$

and $q_4(\ell, u_\alpha)$ denotes the $1-u$ quantile of the random variable $T_4^\ell(\boldsymbol{\varepsilon})$.

The procedure for calculating $q_4(\ell, u_\alpha)$ for $\ell = 2, \dots, \ell_n$ is the same as the procedure for calculating the $q_2(\ell, u_\alpha)$'s.

5.2. *The simulation experiment.* The number of observations $n$ equals 100, $x_i = i/(n+1)$, for $i = 1, \dots, n$, and $\ell_n$ is either equal to 15 or 25.

We consider three distributions of the errors $\varepsilon_i$, with expectation zero and variance 1.

1. The Gaussian distribution: $\varepsilon_i \sim \mathcal{N}(0, 1)$.
2. The type I distribution: $\varepsilon_i$ has density $sf_X(\mu + sx)$, where $f_X(x) = \exp\{-x - \exp(-x)\}$ and where $\mu$ and $s^2$ are the expectation and the variance of a variable $X$ with density $f_X$. This distribution is asymmetrical.
3. The mixture of Gaussian distributions: $\varepsilon_i$ is distributed as $\pi X_1 + (1-\pi)X_2$, where $\pi$ is distributed as a Bernoulli variable with expectation 0.9, $X_1$ and $X_2$ are centered Gaussian variables with variances, respectively, equal to $2.43s$ and $25s$, and $\pi, X_1$ and $X_2$ are independent. The quantity $s$ is chosen such that the variance of $\varepsilon_i$ equals 1. This distribution has heavy tails.

We consider several functions $F$ that are presented below. For each of them, we simulate the observations $Y_i = F(x_i) + \sigma\varepsilon_i$. The values of $\sigma^2$ and of the distance in sup-norm between $F$ and $\mathcal{K}_\nearrow$ are reported in Table 1:

$$d_\infty(F, \mathcal{K}_\nearrow) = \tfrac{1}{2} \sup_{0\le s\le t\le 1} (F(s) - F(t)).$$

Let us comment on the choice of the considered functions.

(a) $F_0(x) = 0$ corresponds to the case for which the quantiles $q(\ell, u_\alpha)$ are calculated.

(b) The function $F_1(x) = 15\mathbf{1}_{x\le 0.5}(x-0.5)^3 + 0.3(x-0.5) - \exp(-250(x-0.25)^2)$ presents a strongly increasing part with a pronounced dip around $x = 1/4$ followed by a nearly flat part on the interval $[1/2, 1]$.

(c) The decreasing linear function $F_2(x) = -ax$, the parameter $a$ being chosen such that $a = 1.5\sigma$.

(d) The function $F_3(x) = -0.2\exp(-50(x-0.5)^2)$ deviates from $F_0$ by a smooth dip while the function $F_4(x) = 0.1\cos(6\pi x)$ deviates from $F_0$ by a cosine function.

(e) The functions $F_5(x) = 0.2x + F_3(x)$ and $F_6(x) = 0.2x + F_4(x)$ deviate from an increasing linear function in the same way as $F_3$ and $F_4$ do from $F_0$.

Let us mention that it is more difficult to detect that $F_5$ (resp. $F_6$) is nonincreasing than to detect that $F_3$ (resp. $F_4$) is. Indeed, adding an increasing function to a function $F$ reduces the distance in sup-norm between $F$ and $\mathcal{K}_{\nearrow}$. This is the reason why the values of $\sigma$ are smaller in the simulation study when we consider the functions $F_5$ and $F_6$.

In Figure 1 we have displayed the functions $F_\ell$ for $\ell = 1, \ldots, 6$ and for each of them one sample simulated with Gaussian errors. The corresponding values of the test statistics $T_{\mathrm{LM}}$ and $T_{\mathrm{LG}}$ for $\alpha = 5\%$ and $\ell_n = 25$ are given. For this simulated sample, it appears that the test based on the statistic $T_{\mathrm{LM}}$ leads to rejection of the null hypothesis in all cases, while the test based on $T_{\mathrm{LG}}$ rejects in all cases except for functions $F_2$ and $F_4$.

The results of the simulation experiment based on 4000 simulations are presented in Tables 2 and 3.

5.3. *Comments on the simulation study.* As expected, the estimated level of the test calculated for the function $F_0(x) = 0$ is (nearly) equal to $\alpha$ when the errors are distributed as Gaussian variables.

When $\ell_n = 25$, the estimated levels of the tests for the mixture and type I distributions are greater than $\alpha$ (see Table 2). Let us recall that when $\ell_n$ is large, we are considering statistics based on the average of the observations on sets $J$ with small cardinality. Therefore, reducing $\ell_n$ improves the robustness to a non-Gaussian error distribution. This is what we get in Table 2 for $\ell_n = 15$. It also appears that the method based on the local means is more robust than the method based on the local gradients, and that both methods are more robust for the type I distribution that is asymmetric but not heavy tailed, than for the mixture distribution.

TABLE 1
*Testing monotonicity: simulated functions $F$, values of $\sigma^2$ and distance in sup-norm between $F$ and $\mathcal{K}_{\nearrow}$*

| $F$ | $\sigma^2$ | $d_\infty(F, \mathcal{K}_{\nearrow})$ |
|---|---|---|
| $F_0(x)$ | 0.01 | 0 |
| $F_1(x)$ | 0.01 | 0.25 |
| $F_2(x)$ | 0.01 | 0.073 |
| $F_3(x)$ | 0.01 | 0.1 |
| $F_4(x)$ | 0.01 | 0.1 |
| $F_5(x)$ | 0.004 | 0.06 |
| $F_6(x)$ | 0.006 | 0.08 |

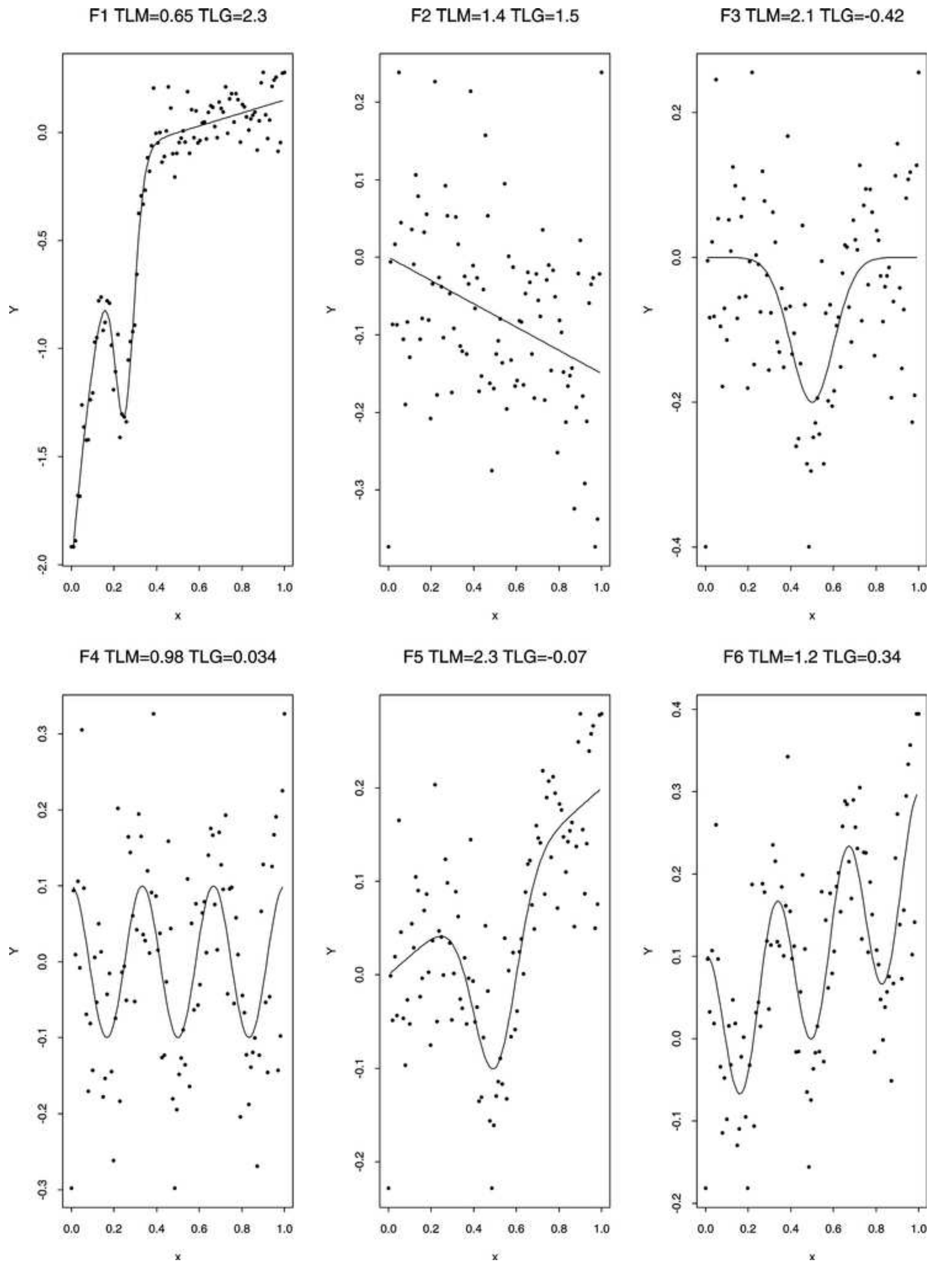

FIG. 1. *For each function $F_\ell$, $\ell = 1, \ldots, 6$, the simulated data $Y_i = F_\ell(x_i) + \sigma\varepsilon_i$ for $i = 1, \ldots, n$ are displayed. The errors $\varepsilon_i$ are Gaussian normalized centered variables. The values of the test statistics $T_{\mathrm{LM}}$ and $T_{\mathrm{LG}}$, with $\alpha = 5\%$, are given for each example.*

Except for the function $F_1$, the estimated power is greater for the procedure based on the local means than for the procedure based on the local gradients (see Table 3). For both procedures the power of the test is larger with $\ell_n = 25$ than with $\ell_n = 15$. However, except for the function $F_1$, the loss of power is less significant for the procedure based on the local means.

5.4. *Comparison with other work.* As expected, the power of our procedure $T_{\mathrm{LG}}$ for the function $F_1$ is similar to that obtained by Hall and Heckman (2000).

The decreasing linear function $F_2(x) = -ax$ has already been studied by Gijbels, Hall, Jones and Koch (2000) with $a = 3\sigma$. They get an estimated power of 77%.

Gijbels, Hall, Jones and Koch (2000) studied the function $0.075F_3/0.2$ with $\sigma = 0.025$ and obtained a simulated power of 98%. With the same function and the same $\sigma$, we get a power equal to 1, for both procedures and for $\ell_n = 15$ and $\ell_n = 25$.

TABLE 2
*Testing monotonicity: levels of the tests based on $T_{\mathrm{LM}}$ and $T_{\mathrm{LG}}$*

| | $\ell_n = 15$ | | $\ell_n = 25$ | |
|---|---|---|---|---|
| **Errors distribution** | $T_{\mathrm{LM}}$ | $T_{\mathrm{LG}}$ | $T_{\mathrm{LM}}$ | $T_{\mathrm{LG}}$ |
| Gaussian | 0.049 | 0.050 | 0.046 | 0.051 |
| Type I | 0.048 | 0.072 | 0.064 | 0.085 |
| Mixture | 0.064 | 0.117 | 0.093 | 0.180 |

TABLE 3
*Testing monotonicity: powers of the tests based on $T_{\mathrm{LM}}$ and $T_{\mathrm{LG}}$ when the errors are Gaussian*

| | $\ell_n = 15$ | | $\ell_n = 25$ | |
|---|---|---|---|---|
| $F$ | $T_{\mathrm{LM}}$ | $T_{\mathrm{LG}}$ | $T_{\mathrm{LM}}$ | $T_{\mathrm{LG}}$ |
| $F_1$ | 0.85 | 0.99 | 0.99 | 1 |
| $F_2$ | 0.96 | 0.96 | 0.99 | 0.99 |
| $F_3$ | 0.99 | 0.73 | 1 | 0.98 |
| $F_4$ | 0.89 | 0.71 | 0.99 | 0.94 |
| $F_5$ | 0.99 | 0.69 | 0.99 | 0.87 |
| $F_6$ | 0.87 | 0.79 | 0.98 | 0.93 |

Gijbels, Hall, Jones and Koch (2000) and Hall and Heckman (2000) calculated the power of their test for the function $F_7(x) = 1 + x - a\exp(-50(x - 0.5)^2)$ for different values of $a$ and $\sigma$. When $a = 0.45$ and $\sigma = 0.05$, we get a power equal to 1 as Gijbels, Hall, Jones and Koch (2000) do. When $a = 0.45$ and $\sigma = 0.1$, we get a power equal to 76% when using the procedure $T_{\mathrm{LM}}$ with $\ell_n = 25$ or $\ell_n = 15$. Gijbels, Hall, Jones and Koch (2000) got 80% and Hall and Heckman (2000) a power larger than 87%.

## 6. Proof of Theorem 1.

*Level of the test.* We first prove that for all $\mathbf{t} \in \mathcal{T}_n$, $q_{\mathbf{t}}(\alpha) > 0$. Indeed, thanks to (33), we have

$$\begin{aligned}
&\mathbb{P}\left[\sqrt{n-d_n}\frac{\langle \boldsymbol{\varepsilon}, \mathbf{t}\rangle}{\|\boldsymbol{\varepsilon} - \Pi_{V_n}\boldsymbol{\varepsilon}\|} - q_{\mathbf{t}}(\alpha) > 0\right] \\
&\qquad \le \mathbb{P}\left[\sup_{\mathbf{t}\in\mathcal{T}_n}\left(\sqrt{n-d_n}\frac{\langle \boldsymbol{\varepsilon}, \mathbf{t}\rangle}{\|\boldsymbol{\varepsilon} - \Pi_{V_n}\boldsymbol{\varepsilon}\|} - q_{\mathbf{t}}(\alpha)\right) > 0\right] \\
&\qquad \le \alpha < \frac{1}{2}.
\end{aligned}$$

Since the random variable $\sqrt{n-d_n}\,\langle \boldsymbol{\varepsilon}, \mathbf{t}\rangle / \|\boldsymbol{\varepsilon} - \Pi_{V_n}\boldsymbol{\varepsilon}\|$ is symmetric (distributed as Student with $n - d_n$ degrees of freedom), we deduce that $q_{\mathbf{t}}(\alpha)$ is positive. In the sequel let us set

$$\hat{\sigma}_n = \|\mathbf{Y} - \Pi_{V_n}\mathbf{Y}\| / \sqrt{n-d_n}.$$

Since for all $\mathbf{f} \in \mathcal{C}$ and $j \in \{1, \ldots, p\}$, $\langle \mathbf{f}, v_j\rangle \le 0$, we have that for all $t \in \mathcal{T}_n$,

$$\langle \mathbf{f}, \mathbf{t}\rangle = \sum_{j=1}^{p} \frac{\lambda_j \langle \mathbf{f}, v_j\rangle}{\|\sum_{j=1}^{p} \lambda_j v_j\|} \le 0.$$

Hence, $\langle \mathbf{Y}, \mathbf{t}\rangle = \langle \mathbf{f}, \mathbf{t}\rangle + \sigma\langle \boldsymbol{\varepsilon}, \mathbf{t}\rangle \le \sigma\langle \boldsymbol{\varepsilon}, \mathbf{t}\rangle$ and therefore for all $\mathbf{f} \in \mathcal{C}$ and $\sigma > 0$,

$$\begin{aligned}
\mathbb{P}_{\mathbf{f},\sigma}[T_\alpha > 0] &\le \mathbb{P}_{\mathbf{f},\sigma}\left[\sup_{\mathbf{t}\in\mathcal{T}_n}\left(\frac{\langle \boldsymbol{\varepsilon}, \mathbf{t}\rangle}{\hat{\sigma}_n/\sigma} - q_{\mathbf{t}}(\alpha)\right) > 0\right] \\
&\le \mathbb{P}_{\mathbf{f},\sigma}\left[\frac{\hat{\sigma}_n}{\sigma} < \sup_{\mathbf{t}\in\mathcal{T}_n}\frac{\langle \boldsymbol{\varepsilon}, \mathbf{t}\rangle}{q_{\mathbf{t}}(\alpha)}\right].
\end{aligned}$$

We now use the following lemma for noncentral $\chi^2$-random variables.

Lemma 1. *For all $u > 0$, $\mathbf{f} \in \mathbb{R}^n$ and $\sigma > 0$*

$$\mathbb{P}_{\mathbf{f},\sigma}[\hat{\sigma}_n < \sigma u] \le \mathbb{P}_{0,1}[\hat{\sigma}_n < u].$$

This lemma states that a noncentral $\chi^2$-random variable is stochastically larger than a $\chi^2$-random variable with the same degrees of freedom. For a proof we refer to Lemma 1 in Baraud, Huet and Laurent (2003a).

Since $\mathcal{T}_n \subset V_n$, the random variables $\langle \boldsymbol{\varepsilon}, \mathbf{t}\rangle$ for $\mathbf{t} \in \mathcal{T}_n$ are independent of $\hat{\sigma}_n$ and thus by conditioning with respect to the $\langle \boldsymbol{\varepsilon}, \mathbf{t}\rangle$'s and using Lemma 1 we get

$$\begin{aligned}\sup_{\sigma>0}\sup_{\mathbf{f}\in\mathcal{C}} \mathbb{P}_{\mathbf{f},\sigma}[T_\alpha > 0] &\le \mathbb{P}_{0,1}\left[\hat{\sigma}_n < \sup_{\mathbf{t}\in\mathcal{T}_n} \frac{\langle \boldsymbol{\varepsilon}, \mathbf{t}\rangle}{q_{\mathbf{t}}(\alpha)}\right]\\ &= \mathbb{P}_{0,1}[T_\alpha > 0] = \alpha.\end{aligned}$$

The reverse inequality being obvious, this concludes the proof of (36).

*Power of the test.* For any $\mathbf{f} \in \mathbb{R}^n$ and $\sigma > 0$

$$\mathbb{P}_{\mathbf{f},\sigma}(T_\alpha \le 0) = \mathbb{P}_{\mathbf{f},\sigma}(\forall\, \mathbf{t} \in \mathcal{T}_n, \langle \mathbf{Y}, \mathbf{t}\rangle \le q_{\mathbf{t}}(\alpha)\hat{\sigma}_n).$$

Setting

$$x_n(\mathbf{f}, \beta) = \frac{\sigma}{\sqrt{n-d_n}} \sqrt{\bar{\chi}^{-1}_{n-d_n,\|\mathbf{f}-\Pi_{V_n}\mathbf{f}\|^2/\sigma^2}(\beta/2)},$$

we have

$$\mathbb{P}_{\mathbf{f},\sigma}(\hat{\sigma}_n > x_n(\mathbf{f}, \beta)) = \beta/2.$$

It follows that for all $\mathbf{f} \in \mathbb{R}^n$ and $\sigma > 0$,

$$\begin{aligned}\mathbb{P}_{\mathbf{f},\sigma}(T_\alpha \le 0) &\le \inf_{t\in\mathcal{T}_n} \mathbb{P}_{\mathbf{f},\sigma}(\langle \mathbf{Y}, \mathbf{t}\rangle \le q_{\mathbf{t}}(\alpha) x_n(\mathbf{f}, \beta)) + \beta/2\\ &\le \inf_{t\in\mathcal{T}_n} \mathbb{P}_{\mathbf{f},\sigma}(\sigma\langle \boldsymbol{\varepsilon}, \mathbf{t}\rangle \le q_{\mathbf{t}}(\alpha) x_n(\mathbf{f}, \beta) - \langle \mathbf{f}, \mathbf{t}\rangle) + \beta/2.\end{aligned}$$

Since $\|t\| = 1$, $\langle \boldsymbol{\varepsilon}, \mathbf{t}\rangle$ is distributed as a standard Gaussian variable, and therefore $\mathbb{P}_{\mathbf{f},\sigma}(T_\alpha \le 0) \le \beta$ as soon as there exists $t \in \mathcal{T}_n$ such that

$$q_{\mathbf{t}}(\alpha) x_n(\mathbf{f}, \beta) - \langle \mathbf{f}, \mathbf{t}\rangle \le -\sigma\bar{\Phi}^{-1}(\beta/2).$$

This concludes the proof of Theorem 1.

**7. Proofs of Propositions 4 and 5.** Let us denote by $\mathcal{I}_r$ the set of increasing sequences of $r+1$ indices in $\{1, \dots, n\}$, that is,

$$\mathcal{I}_r = \{(i_1, \dots, i_{r+1}),\ i_1 < \cdots < i_{r+1},\ i_j \in \{1, \dots, n\}\}. \tag{39}$$

For $\mathbf{i} = (i_1, \dots, i_{r+1}) \in \mathcal{I}_r$ and $\mathbf{v} \in \mathbb{R}^n$ we set

$$\phi_{\mathbf{i}}(\mathbf{v}) = \det\begin{pmatrix} 1 & x_{i_1} & \cdots & x_{i_1}^{r-1} & v_{i_1}\\ 1 & x_{i_2} & \cdots & x_{i_2}^{r-1} & v_{i_2}\\ \vdots & \vdots & \vdots & \vdots & \vdots\\ 1 & x_{i_{r+1}} & \cdots & x_{i_{r+1}}^{r-1} & v_{i_{r+1}}\end{pmatrix}. \tag{40}$$

For $\mathbf{i} = (i, \dots, i+r)$, $\phi_{\mathbf{i}}(\mathbf{v}) = \phi_i(\mathbf{v})$, where $\phi_i(\mathbf{v})$ is defined by (30). For $\mathbf{w}^1, \dots, \mathbf{w}^q$, $q$ vectors of $\mathbb{R}^n$, we set

$$\mathrm{Gram}(\mathbf{w}^1, \dots, \mathbf{w}^q) = \det(G) \qquad \text{where } G = (\langle \mathbf{w}^i, \mathbf{w}^j \rangle)_{1 \le i,j \le q}.$$

Let us define

$$\tilde{\mathcal{C}}_{r,R} = \{\mathbf{f} \in \mathbb{R}^n,\ \forall \mathbf{i} \in \mathcal{I}_r,\ \phi_{\mathbf{i}}(R \star \mathbf{f}\,) \ge 0\}. \tag{41}$$

The proofs of Propositions 4 and 5 rely on the following lemma.

LEMMA 2. *The following equalities hold. First,*

$$\tilde{\mathcal{C}}_{r,R} = \mathcal{C}_{r,R}. \tag{42}$$

*Assume that $\mathbf{f} = (F(x_1), \dots, F(x_n))'$, where $F$ is such that $RF$ is $r$ times differentiable. Then for each $\mathbf{i} \in \mathcal{I}_r$ there exists some $c_{\mathbf{i}} \in ]x_{i_1}, x_{i_{r+1}}[$ such that*

$$\phi_{\mathbf{i}}(R \star \mathbf{f}\,) = \frac{\Lambda(F)(c_{\mathbf{i}})}{r!} \phi_{\mathbf{i}}(\mathbf{x}^r). \tag{43}$$

*For $J \subset \{1, \dots, n\}$ let $\mathbf{t}^*_J$ be defined by (22). We have*

$$-\langle \mathbf{f}, \mathbf{t}^*_J \rangle = N_J^{-1} \sum_{\mathbf{i} \in \mathcal{I}_r \cap J^{r+1}} \phi_{\mathbf{i}}(R \star \mathbf{f}\,) \phi_{\mathbf{i}}(\mathbf{x}^r), \tag{44}$$

*where $N_J = \mathrm{Gram}(\mathbf{1}_J, \mathbf{x}_J, \dots, \mathbf{x}_J^{r-1}) \gamma_J$.*

The proof of the lemma is delayed to the Appendix.

7.1. *Proof of Proposition 4.* The result concerning $\mathcal{K}_{\smile}$ is clear as a function $F$ is nonconcave on $[0,1]$ if and only if for all $x, y, z$ in $[0,1]$ with $x < y < z$ one has

$$\det \begin{pmatrix} 1 & x & F(x) \\ 1 & y & F(y) \\ 1 & z & F(z) \end{pmatrix} \ge 0.$$

Let us now turn to the set $\mathcal{K}_{r,R}$. First note that the $n-r$ linear forms $\mathbf{f} \mapsto \phi_{i,r}(R \star \mathbf{f}\,)$ are independent since the linear space

$$\{\mathbf{f} \in \mathbb{R}^n,\ \forall i \in \{1, \dots, n-r\},\ \phi_{i,r}(R \star \mathbf{f}\,) = 0\},$$

which is generated by

$$\frac{1}{R} \star \mathbf{1}, \frac{1}{R} \star \mathbf{x}, \dots, \frac{1}{R} \star \mathbf{x}^{r-1},$$

is of dimension $r$. Second, the fact that $\mathbf{f}$ belongs to $\tilde{\mathcal{C}}_{r,R}$ is a straightforward consequence of (43) since under the assumption that $F \in \mathcal{K}_{r,R}$, $\Lambda(F)(x) \ge 0$ for all $x$, and since the Vandermonde determinants $\phi_{\mathbf{i}}(\mathbf{x}^r)$ are positive for all $\mathbf{i} \in \mathcal{I}_r$.

7.2. *Proof of Proposition* 5. The result is clear in the case where $\mathcal{C}_{\geq 0}$. For the other cases we use the following lemma.

LEMMA 3. *Let $W$ be the orthogonal complement of the linear space generated by the $\mathbf{v}_j$'s for $j = 1, \dots, p$. If $\mathbf{t}^* \notin W$ satisfies for all $\mathbf{f} \in \mathcal{C}$*

$$\langle \mathbf{t}^* - \Pi_W \mathbf{t}^*, \mathbf{f} \rangle \leq 0,$$

*then*

$$\frac{\mathbf{t}^* - \Pi_W \mathbf{t}^*}{\|\mathbf{t}^* - \Pi_W \mathbf{t}^*\|} \in \mathcal{T}.$$

PROOF. The vector $\mathbf{t}^* - \Pi_W \mathbf{t}^*$ belongs to the linear space generated by the $\mathbf{v}_j$'s and thus one can write $\mathbf{g}^* = \mathbf{t}^* - \Pi_W \mathbf{t}^* = \sum_{j=1}^p \lambda_j \mathbf{v}_j$. It remains to show that the $\lambda_j$'s are nonnegative. Let us fix $j_0 \in \{1, \dots, p\}$ and choose $\mathbf{f}^{j_0}$ in $\mathbb{R}^n$ satisfying $\langle \mathbf{f}^{j_0}, \mathbf{v}_j \rangle = 0$ for all $j \neq j_0$ and $\langle \mathbf{f}^{j_0}, \mathbf{v}_{j_0} \rangle < 0$. Such a vector exists since the $\mathbf{v}_j$'s are linearly independent in $\mathbb{R}^n$. Clearly $\mathbf{f}^{j_0}$ belongs to $\mathcal{C}$ and therefore $\langle \mathbf{f}^{j_0}, \mathbf{g}^* \rangle = \lambda_{j_0} \langle \mathbf{f}^{j_0}, \mathbf{v}_{j_0} \rangle \leq 0$ which constrains $\lambda_{j_0}$ to be nonnegative. This concludes the proof of Lemma 3. □

Let us consider the case where $\mathcal{C} = \mathcal{C}_{\nearrow}$. We apply Lemma 3. In this case $W$ is the linear space generated by $\mathbf{1}$; we get that for all $\ell \in \{2, \dots, \ell_n\}$ and $1 \leq i < j \leq \ell$, $\mathbf{e}_{ij}^\ell$ satisfies $\Pi_W \mathbf{e}_{ij}^\ell = 0$. Moreover $\|\mathbf{e}_{ij}^\ell\| = 1$ and

$$\forall \mathbf{f} \in \mathcal{C}_{\nearrow} \qquad \langle \mathbf{f}, \mathbf{e}_{ij}^\ell \rangle = N_{ij}^\ell (\bar{f}_{J_i^\ell} - \bar{f}_{J_j^\ell}) \leq 0.$$

Let us consider the case where $\mathcal{C} = \mathcal{C}_{\smile}$. In this case, $p = n - 2$ and for all $j = 1, \dots, n-2$,

$$\mathbf{v}_j = (x_{j+1} - x_{j+2})\mathbf{e}_j + (x_{j+2} - x_j)\mathbf{e}_{j+1} + (x_j - x_{j+1})\mathbf{e}_{j+2}.$$

Since $\|\mathbf{e}_{ijk}^\ell\| = 1$, by Lemma 3 it is enough to prove that:

(i) for all $\mathbf{f} \in W$, $\langle \mathbf{f}, \mathbf{e}_{ijk}^\ell \rangle = 0$,
(ii) for all $\mathbf{f} \in \mathcal{C}_{\smile}$, $\langle \mathbf{f}, \mathbf{e}_{ijk}^\ell \rangle \leq 0$.

First note that for all $\mathbf{f} \in \mathbb{R}^n$,

$$\langle \mathbf{f}, \mathbf{e}_{ijk}^\ell \rangle = N_{ijk}^\ell (\bar{f}_{J_j^\ell} - \lambda_{ijk}^\ell \bar{f}_{J_i^\ell} - (1 - \lambda_{ijk}^\ell) \bar{f}_{J_k^\ell}). \tag{45}$$

Clearly if $\mathbf{f} = \mathbf{1}$ or $\mathbf{f} = \mathbf{x}$, $\langle \mathbf{f}, \mathbf{e}_{ijk}^\ell \rangle = 0$ and since by definition of $\mathcal{C}_{\smile}$, $W$ is the linear space generated by $\mathbf{1}$ and $\mathbf{x}$, (i) holds true. Let now $\mathbf{f} \in \mathcal{C}_{\smile}$. There exists some convex function $F$ mapping $[x_1, x_n]$ into $\mathbb{R}$ such that $F(x_i) = f_i$ for all $i = 1, \dots, n$ (take the piecewise linear function verifying this property, e.g.). Let $i < j < k$ and $l \in J_j^\ell$. We set

$$\mu_{ik}^l = \frac{\bar{\mathbf{x}}_{J_k^\ell} - x_l}{\bar{\mathbf{x}}_{J_k^\ell} - \bar{\mathbf{x}}_{J_i^\ell}}.$$

Note that $0 \le \mu_{ik}^l \le 1$ and that

$$x_l = \mu_{ik}^l \bar{\mathbf{x}}_{J_i^\ell} + (1-\mu_{ik}^l)\bar{\mathbf{x}}_{J_k^\ell}.$$

Since $F$ is convex on $[x_1, x_n]$, we have for all $l \in J_j^\ell$,

$$F(x_l) \le \mu_{ik}^l F(\bar{\mathbf{x}}_{J_i^\ell}) + (1-\mu_{ik}^l)F(\bar{\mathbf{x}}_{J_k^\ell}) \le \mu_{ik}^l \bar{\mathbf{f}}_{J_i^\ell} + (1-\mu_{ik}^l)\bar{\mathbf{f}}_{J_k^\ell}.$$

Note that $\sum_{l\in J_j^\ell} \mu_{ik}^l/|J_j^\ell| = \lambda_{ijk}^\ell$. We derive from the above inequality that

$$\bar{f}_{J_j^\ell} = \frac{1}{|J_j^\ell|}\sum_{l\in J_j^\ell} F(x_l) \le \lambda_{ijk}^\ell \bar{f}_{J_i^\ell} + (1-\lambda_{ijk}^\ell)\bar{f}_{J_k^\ell},$$

which, thanks to (45), leads to (ii).

Let us consider the case where $\mathcal{C} = \mathcal{C}_{r,R}$. By Lemma 2 we know that $\tilde{\mathcal{C}}_{r,R} = \mathcal{C}_{r,R}$ and therefore for each $\mathbf{i} \in \mathcal{I}_r$, the linear form $\mathbf{f} \mapsto \phi_{\mathbf{i}}(R\star\mathbf{f})$ is a linear combination of the linear forms $\mathbf{f} \mapsto \phi_i(R\star\mathbf{f})$ with $i = 1, \dots, n-r$. Consequently, if $\mathbf{w} \in W$, then for all $\mathbf{i} \in \mathcal{I}_r$, $\phi_{\mathbf{i}}(R\star\mathbf{w}) = 0$. For each $J \subset \{1,\dots,n\}$, $\mathbf{t}_J^*$ defined by (22) satisfies $\|\mathbf{t}_J^*\| = 1$. By applying (44) with $\mathbf{f} = \mathbf{w}$, we get $\langle \mathbf{w}, \mathbf{t}_J^*\rangle \le 0$ for all $\mathbf{w} \in \mathcal{C}_{r,R}$ and $\langle \mathbf{w}, \mathbf{t}_J^*\rangle = 0$ for all $\mathbf{w} \in W$. Consequently, by Lemma 3, $\mathbf{t}_J^*$ belongs to $\mathcal{T}$.

## 8. Proof of Proposition 1.

8.1. *Proof for* $(T_\alpha, \mathcal{C}) = (T_{\alpha,1}, \mathcal{C}_{\ge 0})$. We prove the proposition by applying Theorem 1. We decompose the proof into six steps.

STEP 1. *For all integer* $N \ge 1$, *let* $\bar{T}_N^{-1}(u)$ *denote the* $1-u$ *quantile of a Student random variable with* $N$ *degrees of freedom. We have for all* $u \in ]0,1[$,

$$\bar{T}_N^{-1}(u) \le 1 + C\left\{\log^{1/4}\left(\frac{1}{u}\right) + \log^{1/2}\left(\frac{1}{u}\right)\exp\left(\frac{2}{N}\log\left(\frac{1}{u}\right)\right)\right\} \tag{46}$$

*for some absolute constant* $C > 0$.

PROOF. Let $\bar{F}_{1,N}^{-1}(u)$ denote the $1-u$ quantile of a Fisher variable with one and $N$ degrees of freedom. Then

$$\bar{T}_N^{-1}(u) = \sqrt{\bar{F}_{1,N}^{-1}(u)}.$$

It follows from Lemma 1 in Baraud, Huet and Laurent (2003a) that for all $u \in ]0,1[$, $N \ge 1$,

$$\bar{F}_{1,N}^{-1}(u) \le 1 + 2\sqrt{2}\log^{1/2}\left(\frac{1}{u}\right) + \frac{3N}{2}\left\{\exp\left(\frac{4}{N}\log\left(\frac{1}{u}\right)\right) - 1\right\}.$$

Using the inequality $\exp(x)-1\le x\exp(x)$ which holds for all $x>0$, we obtain

$$\bar F_{1,N}^{-1}(u)\le 1+2\sqrt{2}\log^{1/2}\left(\frac{1}{u}\right)+6\log\left(\frac{1}{u}\right)\exp\left(\frac{4}{N}\log\left(\frac{1}{u}\right)\right),$$

and since $\sqrt{a+b}\le\sqrt{a}+\sqrt{b}$ for all $a>0$ and $b>0$,

$$\sqrt{\bar F_{1,N}^{-1}(u)}\le 1+C\left\{\log^{1/4}\left(\frac{1}{u}\right)+\log^{1/2}\left(\frac{1}{u}\right)\exp\left(\frac{2}{N}\log\left(\frac{1}{u}\right)\right)\right\}$$

for some absolute constant $C>0$. □

STEP 2. *For all* $\ell\in\{1,\dots,\ell_n\}$, $t\in\mathcal{T}_{n,1}^{\ell}$, *we have*

$$q_{\mathbf{t}}(\alpha)=q_1(\ell,u_\alpha)\le C(\alpha)\sqrt{\log(n)}. \tag{47}$$

PROOF. On the one hand, by definition of $q_1(\ell,\cdot)$,

$$\alpha=\mathbb{P}_{0,1}(T_{\alpha,1}>0)\le\sum_{\ell=1}^{\ell_n}\mathbb{P}(T_1^{\ell}(\boldsymbol{\varepsilon})-q_1(\ell,u_\alpha)>0)\le\ell_n u_\alpha,$$

and thus

$$u_\alpha\ge\alpha/\ell_n. \tag{48}$$

On the other hand, for all $\ell\in\{1,\dots,\ell_n\}$ and $J\in\mathcal{J}^{\ell}$, the random variables

$$U_J=\frac{-\sum_{i\in J}\varepsilon_i}{\|\boldsymbol{\varepsilon}-\Pi_{V_n}\boldsymbol{\varepsilon}\|}\sqrt{\frac{n-d_n}{|J|}}$$

being distributed as Student variables with $n-d_n$ degrees of freedom, we have that

$$\mathbb{P}\left(T_1^{\ell}(\boldsymbol{\varepsilon})>\bar T_{n-d_n}^{-1}\left(\frac{u_\alpha}{|\mathcal{J}^{\ell}|}\right)\right)\le\sum_{J\in\mathcal{J}^{\ell}}\mathbb{P}\left(U_J>\bar T_{n-d_n}^{-1}\left(\frac{u_\alpha}{|\mathcal{J}^{\ell}|}\right)\right)\le u_\alpha \tag{49}$$

and thus $q_1(\ell,u_\alpha)\le\bar T_{n-d_n}^{-1}(u_\alpha/|\mathcal{J}^{\ell}|)$. This inequality together with (48) and (46) leads to (47), as $|\mathcal{J}^{\ell}|\le\ell_n\le n/2$ and $n-d_n=n-\ell_n\ge n/2$. □

STEP 3. *For all* $\mathbf{f}=(F(x_1),\dots,F(x_n))'$ *with* $F\in\mathcal{H}_s(L)$,

$$\frac{\|\mathbf{f}-\Pi_{V_{n,\mathrm{cste}}}\mathbf{f}\|^2}{n}\le C(s)L^2n^{-2s}. \tag{50}$$

PROOF. Note that the vector

$$\tilde{\mathbf{f}}=\sum_{k=1}^{\ell_n}F(\bar x_{J_k})\mathbf{1}_{J_k}$$

belongs to $V_{n,\mathrm{cste}}$ and therefore

$$\begin{aligned}\|\mathbf{f}-\Pi_{V_{n,\mathrm{cste}}}\mathbf{f}\|^2 &\le \|\mathbf{f}-\tilde{\mathbf{f}}\|^2 \\ &= \sum_{k=1}^{\ell_n}\sum_{i\in J_k}(F(x_i)-F(\bar{x}_{J_k}))^2 \\ &\le \sum_{k=1}^{\ell_n}\sum_{i\in J_k}L^2\ell_n^{-2s} \\ &= nL^2\ell_n^{-2s}.\end{aligned}$$

Noting that $\ell_n = d_n \ge n/4$, we get (50). □

STEP 4. *Assuming that $n \ge (L/\sigma)^{1/s}$, there exists some constant $C$ depending on $s$ and $\beta$ only such that*

$$\frac{\bar{\chi}^{-1}_{n-d_n,\|\mathbf{f}-\Pi_{V_{n,\mathrm{cste}}}\mathbf{f}\|^2/\sigma^2}(\beta/2)}{n-d_n} \le C. \tag{51}$$

PROOF. Using the inequality due to Birgé (2001) on the quantiles of noncentral $\chi^2$, we have that

$$\bar{\chi}^{-1}_{n-d_n,a^2}(\beta/2) \le n-d_n+a^2+2\sqrt{(n-d_n+2a^2)\log(2/\beta)}+2\log(2/\beta).$$

Setting $a = \|\mathbf{f}-\Pi_{V_{n,\mathrm{cste}}}\mathbf{f}\|/\sigma$ and using (50), we derive that

$$\bar{\chi}^{-1}_{n-d_n,a^2}(\beta/2)/(n-d_n) \le C(\beta,s). \tag{52}$$ □

STEP 5. *Under the assumption of Step 4, for all $\mathbf{t} \in \mathcal{T}_n$,*

$$v_{\mathbf{t}}(\mathbf{f},\beta) \le \kappa^*\sqrt{\log(n)}\sigma,$$

*for some constant $\kappa^*$ depending on $\alpha,\beta$ and $s$ only.*

PROOF. We recall that

$$v_{\mathbf{t}}(\mathbf{f},\beta) = \left(q_{\mathbf{t}}(\alpha)\frac{1}{\sqrt{n-d_n}}\sqrt{\bar{\chi}^{-1}_{n-d_n,\|\mathbf{f}-\Pi_{V_n}\mathbf{f}\|^2/\sigma^2}(\beta/2)}+\bar{\Phi}^{-1}(\beta/2)\right)\sigma.$$

We conclude by using the elementary inequality

$$\bar{\Phi}^{-1}(\beta/2) \le \sqrt{2\log(2/\beta)},$$

and by gathering (47) and (51). □

We conclude the proof with this final step.

STEP 6. *There exists a constant $\kappa$ depending on $\alpha, \beta$ and $s$ only, such that if $n$ is large enough and $F$ satisfies*

$$\min_{x\in[0,1]} F(x) \le -\kappa\rho_n, \tag{53}$$

*then there exists $\mathbf{t}^* \in \mathcal{T}_n$ such that*

$$\langle \mathbf{f}, \mathbf{t}^*\rangle \ge v_{\mathbf{t}^*}(\mathbf{f}, \beta). \tag{54}$$

PROOF. Since $F \in \mathcal{H}_s(L)$, under Assumption (53) there exists $j \in \{1, 2, \ldots, n\}$ such that

$$F(j/n) \le -\kappa\rho_n + Ln^{-s}.$$

For $n$ large enough, $Ln^{-s} \le \kappa\rho_n/2$, hence $F(j/n) \le -\kappa\rho_n/2$.

Let us take $\kappa$ satisfying

$$\frac{\kappa}{4} = (2\kappa^*)^{2s/(1+2s)},$$

where $\kappa^*$ is defined at Step 5.

Let us define

$$\ell(n) = \left[\left(\frac{4L}{\kappa\rho_n}\right)^{1/s}\right], \tag{55}$$

and $J$ as the element of $\mathcal{J}^{\ell(n)}$ containing $j$. Note that for $n$ large enough, $\ell(n) \in \{1, \ldots, \ell_n\}$.

Now, for all $k \in J$, since $F \in \mathcal{H}_s(L)$,

$$\begin{aligned} f_k = F(x_k) &= -F(x_j) + F(x_j) + F(x_k) \\ &\le -\kappa\rho_n/2 + L|x_k - x_j|^s \\ &\le -\kappa\rho_n/2 + L\ell(n)^{-s} \\ &\le -\kappa\rho_n/4 \end{aligned}$$

and thus, by taking $\mathbf{t}^* \in \mathcal{T}_{n,1}$ as

$$\mathbf{t}^* = -\frac{1}{\sqrt{|J|}} \sum_{i\in J} \mathbf{e}_i,$$

we derive that

$$\begin{aligned} \langle \mathbf{f}, \mathbf{t}^*\rangle &= -\sqrt{|J|}\bar{f}_J \\ &\ge \sqrt{|J|}\kappa\rho_n/4. \end{aligned}$$

By construction of the partition of the data, we have for all positive integers $p \le q \le r$ that

$$\left[\frac{r}{q}\right] \le |I^r_{p,q}| \le \left[\frac{r}{q}\right] + 1. \tag{56}$$

For all $j \in \{1, \dots, \ell(n)\}$, $J = J_j^{\ell(n)}$ [see (8)] is a union of $|I_{j,\ell(n)}^{\ell_n}| \geq [\ell_n/\ell(n)]$ disjoint sets of cardinality at least $[n/\ell_n]$. Hence

$$|J_j^{\ell(n)}| \geq \left[\frac{n}{\ell_n}\right]\left[\frac{\ell_n}{\ell(n)}\right] \geq \frac{n}{4\ell(n)}$$

since $[x] \geq x/2$ for all $x \geq 1$. Therefore we get

$$|J| \geq \frac{n}{4\ell(n)} \geq \frac{n}{4}\left(\kappa \frac{\rho_n}{4L}\right)^{1/s} \tag{57}$$

using (55).

This implies that

$$\langle \mathbf{f}, \mathbf{t}^* \rangle \geq \frac{\sqrt{n}}{8}\left(\frac{\kappa\rho_n}{4L}\right)^{1/(2s)} \kappa\rho_n \geq \kappa^* \sigma \sqrt{\log(n)}$$

by definition of $\kappa$. $\square$

8.2. *Proof for* $(T_\alpha, \mathcal{C}) = (T_{\alpha,2}, \mathcal{C}_\nearrow)$. We follow the proof of Theorem 1 for $(T_\alpha, \mathcal{C}) = (T_{\alpha,1}, \mathcal{C}_{\geq 0})$: the results of Steps 1–5 still hold. The proof of Step 2 differs in the following way: (49) becomes

$$\mathbb{P}\left(T_2^\ell(\boldsymbol{\varepsilon}) > \bar{T}_{n-d_n}^{-1}\left(\frac{u_\alpha}{|\mathcal{T}_{n,2}^\ell|}\right)\right)$$
$$\leq \sum_{1 \leq i < j \leq \ell} \mathbb{P}\left(\frac{\langle \boldsymbol{\varepsilon}, \mathbf{e}_{ij}^\ell \rangle}{\|\boldsymbol{\varepsilon} - \Pi_{V_{n,\mathrm{cste}}}\boldsymbol{\varepsilon}\|/(n-\ell_n)} > \bar{T}_{n-d_n}^{-1}\left(\frac{u_\alpha}{|\mathcal{T}_{n,2}^\ell|}\right)\right) \leq u_\alpha.$$

We conclude the proof of Step 2 by noticing that for all $\ell \in \{1, \dots, \ell_n\}$, $|\mathcal{T}_{n,2}^\ell|$ is bounded from above by $n^2/4$.

Step 6. *For $n$ large enough, under the assumption that*

$$\inf_{G \in \mathcal{K}_\nearrow} \|F - G\|_\infty \geq \kappa\rho_n, \tag{58}$$

*there exists $\mathbf{t}^* \in \mathcal{T}_{n,2}$, such that $\langle \mathbf{t}^*, f \rangle \geq v_{\mathbf{t}^*}(\mathbf{f}, \beta)$.*

Proof. Let us first remark that

$$\inf_{G \in \mathcal{K}_\nearrow} \|F - G\|_\infty \leq \sup_{0 \leq x \leq y \leq 1} (F(x) - F(y)).$$

Indeed, let $G^* \in \mathcal{K}_\nearrow$ be defined as

$$G^*(y) = \sup_{0 \leq x \leq y} F(x).$$

Then

$$\inf_{G\in\mathcal{K}_{\nearrow}} \|F-G\|_\infty \le \|F-G^*\|_\infty = \sup_{0\le x\le y\le 1}(F(x)-F(y)).$$

Hence, under (58), there exists $x<y$ such that $F(x)-F(y)\ge \kappa\rho_n$. Since $F\in\mathcal{H}_s(L)$, if $|x_i-x|\le 1/n$ and $|x_j-y|\le 1/n$, then

$$F(x_i)-F(x_j)\ge \kappa\rho_n - 2Ln^{-s}\ge \kappa\rho_n/2$$

for $n$ large enough. Hence, there exists $1\le i<j\le n$ such that $F(x_i)-F(x_j)\ge \kappa\rho_n/2$.

Let us set

$$\ell(n)=\left[\left(\frac{8L}{\kappa\rho_n}\right)^{1/s}\right],$$

which belongs to $\{1,\dots,\ell_n\}$ at least for $n$ large enough. Let $I$ and $J$ be the elements of $\mathcal{J}^{\ell(n)}$ satisfying $i\in I$ and $j\in J$.

Arguing as in Step 6 of Section 8.1, since $F\in\mathcal{H}_s(L)$,

$$\bar f_I \ge F(x_i)-L\ell(n)^{-s}\quad\text{and}\quad \bar f_J\le F(x_j)+L\ell(n)^{-s}$$

and we deduce that

$$\bar f_I-\bar f_J\ge \kappa\rho_n/2-2L\ell(n)^{-s}\ge \kappa\rho_n/4.$$

This implies that there exists $1\le i^*<j^*\le \ell(n)$ with $I=J_{i^*}^{\ell(n)}$ and $J=J_{j^*}^{\ell(n)}$, such that

$$\langle \mathbf{e}_{i^*j^*}^{\ell(n)}, f\rangle = N_{i^*j^*}^{\ell(n)}(\bar f_I-\bar f_J)\ge N_{i^*j^*}^{\ell(n)}\frac{\kappa\rho_n}{4}.$$

Using (56), and since $\ell_n=[n/2]$, we have that for all $K\in\mathcal{J}^{\ell(n)}$,

$$2\left[\frac{\ell_n}{\ell(n)}\right]\le |K|\le 3\left(\left[\frac{\ell_n}{\ell(n)}\right]+1\right],$$

which implies that

$$N_{i^*j^*}^{\ell(n)}=\sqrt{\frac{|I||J|}{|I|+|J|}}\ge C\sqrt{\frac{\ell_n}{\ell(n)}}.$$

We now conclude as in the proof of Step 6 by taking $\mathbf{t}^*=\mathbf{e}_{i^*j^*}^{\ell(n)}$. $\square$

8.3. *Proof of Theorem 1 for* $(T_\alpha, \mathcal{C}) = (T_{\alpha,3}, \mathcal{C}_\smile)$. We follow the proof of Theorem 1 for the case $(T_\alpha, \mathcal{C}) = (T_{\alpha,1}, \mathcal{C}_{\geq 0})$: the results of Steps 1–5 still hold. Nevertheless, the proof of Step 2 differs in the following way: (49) becomes

$$\mathbb{P}\left(T_3^\ell(\boldsymbol{\varepsilon}) > \bar{T}_{n-d_n}^{-1}\left(\frac{u_\alpha}{|\mathcal{T}_{n,3}^\ell|}\right)\right)$$
$$\leq \sum_{1\leq i<j<k\leq \ell} \mathbb{P}\left(\frac{\langle \boldsymbol{\varepsilon}, \mathbf{e}_{ijk}^\ell\rangle}{\|\boldsymbol{\varepsilon} - \Pi_{V_{n,\mathrm{cste}}}\boldsymbol{\varepsilon}\|/(n-\ell_n)} > \bar{T}_{n-d_n}^{-1}\left(\frac{u_\alpha}{|\mathcal{T}_{n,3}^\ell|}\right)\right) \leq u_\alpha.$$

We conclude the proof of Step 2 by noticing that for all $\ell \in \{1, \dots, \ell_n\}$, $|\mathcal{T}_{n,3}^\ell|$ is bounded from above by $n^3/8$.

STEP 6. *For $n$ large enough, under the assumption that*

$$\inf_{G\in\mathcal{K}_\smile} \|F - G\|_\infty \geq \kappa\rho_n, \tag{59}$$

*there exists $\mathbf{t}^* \in \mathcal{T}_{n,3}$ such that*

$$\langle \mathbf{t}^*, f\rangle \geq v_{\mathbf{t}^*}(\mathbf{f}, \beta).$$

PROOF. We decompose the proof into three parts.

PART 1. *For $n$ large enough, and all $F \in \mathcal{H}_s(L)$ satisfying (59), we have*

$$\inf_{\mathbf{g}\in\mathcal{C}_\smile} \|\mathbf{f} - \mathbf{g}\|_\infty \geq \kappa\rho_n/4,$$

*with $\mathbf{f} = (F(x_1), \dots, F(x_n))'$.*

PROOF. We first prove the following inequality:

$$\inf_{G\in\mathcal{K}_\smile} \|F - G\|_\infty \leq 2Ln^{-s} + 3 \inf_{\mathbf{g}\in\mathcal{C}_\smile} \|\mathbf{f} - \mathbf{g}\|_\infty. \tag{60}$$

Part 1 derives obviously from this inequality.

For all $\mathbf{g} \in \mathcal{C}_\smile$, we consider the function $G_{\mathbf{g}} \in \mathcal{K}_\smile$ defined as the piecewise linear function such that for all $i$, $G_{\mathbf{g}}(x_i) = g_i$ and such that $G_{\mathbf{g}}$ is affine on the interval $[0, x_2]$. Then $\inf_{G\in\mathcal{K}_\smile} \|F - G\|_\infty \leq \|F - G_{\mathbf{g}}\|_\infty$. Moreover, by setting $x_0 = 0$ and $g_0 = G_{\mathbf{g}}(0)$,

$$\begin{aligned}
&\|F - G_{\mathbf{g}}\|_\infty \\
&\quad= \sup_{i\in\{1,\dots,n\}} \sup_{x\in[x_{i-1},x_i]} |F(x) - G_{\mathbf{g}}(x)| \\
&\quad\leq \sup_{i\in\{1,\dots,n\}} \sup_{x\in[x_{i-1},x_i]} |F(x) - F(x_i) + F(x_i) - G_{\mathbf{g}}(x_i) + G_{\mathbf{g}}(x_i) - G_{\mathbf{g}}(x)| \\
&\quad\leq Ln^{-s} + \|f - g\|_\infty + \sup_{i\in\{1,\dots,n\}} |g_{i-1} - g_i|,
\end{aligned}$$

since $\sup_{x\in[x_{i-1},x_i]} |G_{\mathbf{g}}(x_i) - G_{\mathbf{g}}(x)| = |G_{\mathbf{g}}(x_i) - G_{\mathbf{g}}(x_{i-1})|$ ($G$ is linear on $[x_{i-1}, x_i]$). In addition, noticing that $|g_1 - g_0| = |g_2 - g_1|$,

$$\begin{aligned}\sup_{i\in\{1,\dots,n\}} |g_i - g_{i-1}| &\le \sup_{i\in\{2,\dots,n\}} |g_i - f_i + f_i - f_{i-1} + f_{i-1} - g_{i-1}| \\ &\le 2\|\mathbf{f} - \mathbf{g}\|_\infty + Ln^{-s}.\end{aligned}$$

This concludes the proof of (60). □

PART 2. *For all* $\mathbf{f} \in \mathbb{R}^n$,

$$\inf_{\mathbf{g}\in\mathcal{C}_\smile} \|\mathbf{f} - \mathbf{g}\|_\infty \le \max_{1\le i<j<k\le n} \left( f_j - \frac{x_k - x_j}{x_k - x_i} f_i - \frac{x_j - x_i}{x_k - x_i} f_k \right)_+, \tag{61}$$

*where for* $x \in \mathbb{R}$, $(x)_+ = x\mathbf{1}_{x>0}$ *denotes the positive part of* $x$.

PROOF. Let us define $\mathbf{g}^* \in \mathcal{C}_\smile$ as follows: $g_1^* = f_1$ and for $i = 1, \dots, n-1$,

$$g_{i+1}^* = g_i^* + \inf\left\{ \frac{f_k - g_i^*}{x_k - x_i},\ k > i \right\} (x_{i+1} - x_i).$$

In words, if $F_{\text{lin}}$ denotes the piecewise linear function on $[x_1, x_n]$ taking the value $f_i$ at $x_i$, then $\mathbf{g}^*$ is the vector $(G_{\text{lin}}^*(x_1), \dots, G_{\text{lin}}^*(x_n))'$, where $G_{\text{lin}}^*$ is the largest convex function satisfying for all $u \in [x_1, x_n]$ $G_{\text{lin}}^*(u) \le F_{\text{lin}}(u)$. Note that the function $G_{\text{lin}}^*$ is also piecewise linear and satisfies that for all $j \in \{1, \dots, n\}$ such that $F_{\text{lin}}(x_j) - G_{\text{lin}}^*(x_j) > 0$, there exist $1 \le i < j < k \le n$ such that

$$F_{\text{lin}}(x_j) - G_{\text{lin}}^*(x_j) = f_j - \frac{x_k - x_j}{x_k - x_i} f_i - \frac{x_j - x_i}{x_k - x_i} f_k.$$

Consequently,

$$\begin{aligned}\|\mathbf{f} - \mathbf{g}^*\|_\infty &= \max_{j=1,\dots,n} (F_{\text{lin}}(x_j) - G_{\text{lin}}^*(x_j)) \\ &\le \max_{1\le i<j<k\le n} \left( f_j - \frac{x_k - x_j}{x_k - x_i} f_i - \frac{x_j - x_i}{x_k - x_i} f_k \right)_+.\end{aligned}$$ □

PART 3. *Let* $\kappa' = \kappa/4$. *We set*

$$\ell(n) = 1 + \left[ \left( \frac{6L}{\kappa' \rho_n} \right)^{1/s} \right].$$

*If there exist* $1 \le i < j < k \le n$ *such that*

$$f_j - \frac{x_k - x_j}{x_k - x_i} f_i - \frac{x_j - x_i}{x_k - x_i} f_k \ge \kappa' \rho_n,$$

*then there exist* $I = J_{i^*}^{\ell(n)}$, $J = J_{j^*}^{\ell(n)}$ *and* $K = J_{k^*}^{\ell(n)}$ *with* $i^* < j^* < k^*$, *such that*

$$\bar f_J - \frac{\bar x_K - \bar x_J}{\bar x_K - \bar x_I} \bar f_I - \frac{\bar x_J - \bar x_I}{\bar x_K - \bar x_I} \bar f_K \ge \kappa' \rho_n / 4. \tag{62}$$

PROOF. Note that

$$\ell(n) \geq \left(\frac{6L}{\kappa'\rho_n}\right)^{1/s} \tag{63}$$

and that for $n$ large enough, $\ell(n) \in \{1, \dots, \ell_n\}$.

In the sequel, we shall use the following inequalities:

$$\begin{aligned} \forall E \in \{I, J, K\} \qquad & \max_{l,l' \in E} |x_l - x_{l'}| \leq \frac{1}{\ell(n)} \quad \text{and} \\ & \max_{l \in E} |x_l - \bar{x}_E| \leq \frac{1}{2\ell(n)}, \end{aligned} \tag{64}$$

and the following notation:

$$\lambda = \frac{x_k - x_j}{x_k - x_i}, \qquad \bar{\lambda} = \frac{\bar{x}_K - \bar{x}_J}{\bar{x}_K - \bar{x}_I}, \qquad \Delta = \bar{f}_J - \bar{\lambda}\bar{f}_I - (1 - \bar{\lambda})\bar{f}_K.$$

We bound $\Delta$ from below as follows:

$$\begin{aligned} \Delta &= f_j - \lambda f_i - (1 - \lambda) f_k \\ &\quad + \bar{f}_J - f_j + \lambda f_i - \bar{\lambda}\bar{f}_I + (1 - \lambda) f_k - (1 - \bar{\lambda})\bar{f}_K \\ &\geq \kappa\rho_n + \bar{f}_J - f_j + (\lambda - \bar{\lambda}) f_i - \bar{\lambda}(\bar{f}_I - f_i) + (\bar{\lambda} - \lambda) f_k - (1 - \bar{\lambda})(\bar{f}_K - f_k) \\ &\geq \kappa\rho_n - 2\max\{|\bar{f}_I - f_i|, |\bar{f}_J - f_j|, |\bar{f}_K - f_k|\} - |\lambda - \bar{\lambda}||f_i - f_k|. \end{aligned}$$

Let us now bound from above the quantities

$$|f_i - f_k|, \qquad \max\{|\bar{f}_I - f_i|, |\bar{f}_J - f_j|, |\bar{f}_K - f_k|\}, \qquad |\lambda - \bar{\lambda}|.$$

Since $F \in \mathcal{H}_s(L)$, we have that

$$|f_i - f_k| = |F(x_i) - F(x_k)| \leq L|x_k - x_i|^s, \tag{65}$$

and by using (64) that

$$\max\{|\bar{f}_I - f_i|, |\bar{f}_J - f_j|, |\bar{f}_K - f_k|\} \leq L\ell(n)^{-s}. \tag{66}$$

For each $(l, E) \in \{(i, I), (j, J), (k, K)\}$, let

$$h_l = \bar{x}_E - x_l.$$

We have

$$\begin{aligned} \bar{\lambda} &= \frac{x_k - x_j + h_k - h_j}{x_k - x_i + h_k - h_i} \\ &= \lambda\left(\frac{1 + (h_k - h_j)/(x_k - x_j)}{1 + (h_k - h_i)/(x_k - x_i)}\right) \\ &= \lambda\left(1 + \frac{(h_k - h_j)/(x_k - x_j) - (h_k - h_i)/(x_k - x_i)}{1 + (h_k - h_i)/(x_k - x_i)}\right), \end{aligned}$$

and as from (64) $\max\{|h_k - h_j|, |h_k - h_i|\} \le \ell(n)^{-1}$, we deduce that

$$\begin{aligned} |\bar\lambda - \lambda| &= |\lambda| \left| \frac{(h_k - h_j)/(x_k - x_i) - (h_k - h_i)/(x_k - x_i)}{1 + (h_k - h_i)/(x_k - x_i)} \right| \\ &\le \frac{2\delta}{|1-\delta|}, \end{aligned} \tag{67}$$

where

$$\delta = \frac{1}{\ell(n)|x_k - x_i|}. \tag{68}$$

In order to bound $\delta$ from above, note that since $F \in \mathcal{H}_s(L)$,

$$\begin{aligned} \kappa' \rho_n &\le f_j - \lambda f_i - (1-\lambda) f_k \\ &= \lambda(F(x_j) - F(x_i)) + (1-\lambda)(F(x_j) - F(x_k)) \\ &\le L \max\{|x_j - x_i|^s, |x_k - x_j|^s\} \end{aligned}$$

and therefore

$$\begin{aligned} |x_k - x_i| &\ge \max\{|x_j - x_i|, |x_k - x_j|\} \\ &= \{\max\{|x_j - x_i|^s, |x_k - x_j|^s\}\}^{1/s} \\ &\ge \left(\frac{\kappa' \rho_n}{L}\right)^{1/s}. \end{aligned}$$

Thus, we deduce by (63) and the fact that $s \in ]\,0,1]$ that

$$\delta \le \frac{L^{1/s}}{(\kappa' \rho_n)^{1/s} \ell(n)} \le \frac{1}{6}. \tag{69}$$

By gathering (65)–(67), we get

$$\Delta \ge \kappa' \rho_n - 2L\ell(n)^{-s} - 2L \frac{\delta}{1-\delta} |x_k - x_i|^s.$$

By using (68), (69) and (63) we finally get

$$\begin{aligned} \Delta &= \kappa' \rho_n - 2L\ell(n)^{-s} - 2L\ell(n)^{-s} \frac{\delta^{1-s}}{1-\delta} \\ &\ge \kappa' \rho_n \left\{ 1 - \frac{1}{3}\left(1 + \frac{1}{1 - 1/6}\right)\right\} \\ &\ge \kappa' \rho_n / 4. \end{aligned}$$

□

Let us now conclude the proof of Step 6. Under the assumption that

$$\inf_{\mathbf{g} \in \mathcal{C}_{\smile}} \|\mathbf{f} - \mathbf{g}\|_\infty \ge \kappa' \rho_n$$

we know from (61) that there exists $i < j < k$ such that

$$f_j - \frac{x_k - x_j}{x_k - x_i} f_i - \frac{x_j - x_i}{x_k - x_i} f_k \geq \kappa' \rho_n,$$

and from (62) that there exist $I = J_{i^*}^{\ell(n)}$, $J = J_{j^*}^{\ell(n)}$ and $K = J_{k^*}^{\ell(n)}$ with $i^* < j^* < k^*$ such that

$$\langle \mathbf{f}, \mathbf{e}_{i^*j^*k^*}^{\ell(n)} \rangle = N_{i^*j^*k^*}^{\ell(n)} (\bar{f}_J - \lambda_{i^*j^*k^*}^{\ell(n)} \bar{f}_I - (1 - \lambda_{i^*j^*k^*}^{\ell(n)}) \bar{f}_K) \geq \frac{N_{i^*j^*k^*}^{\ell(n)} \kappa' \rho_n}{4}.$$

Noting that for all $E \in \{I, J, K\}$

$$|E| \geq 2 \left[ \frac{\ell_n}{\ell(n)} \right] \geq \frac{\ell_n}{\ell(n)} \geq \frac{n}{4\ell(n)},$$

and that $\|\mathbf{e}_{i^*j^*k^*}^{\ell(n)}\|^2 \leq 1/|I| + 1/|J| + 1/|K|$, we have that

$$N_{i^*j^*k^*}^{\ell(n)} \geq \sqrt{\frac{1}{|I|^{-1} + |J|^{-1} + |K|^{-1}}} \geq \sqrt{\frac{n}{12\ell(n)}}.$$

As $\ell(n) \leq 2(12L/(\kappa' \rho_n))^{1/s}$ at least for $n$ large enough, we deduce that

$$N_{i^*j^*k^*}^{\ell(n)} \geq \sqrt{\frac{n(\kappa' \rho_n)^{1/s}}{8(12L)^{1/s}}}.$$

Consequently, we get

$$\begin{aligned} \langle \mathbf{f}, \mathbf{e}_{i^*j^*k^*}^{\ell(n)} \rangle &\geq \sqrt{(\kappa' \rho_n)^{(1+2s)/s} \frac{n}{12^{1/s} 128 L^{1/s}}} \\ &\geq \kappa^* \sqrt{\log(n)}\, \sigma, \end{aligned}$$

for $\kappa'$ suitably chosen. It remains to take $\mathbf{t}^* = \mathbf{e}_{i^*j^*k^*}^{\ell(n)} \in \mathcal{T}_{n,3}$ to complete the proof.

□

**9. Proof of Proposition 3.** The proof of Proposition 3 is divided into two parts. In Section 9.1 we show that if (24) or (25) holds, then $\mathbb{P}_{F,\sigma}(T_\alpha > 0) \geq 1 - \beta$. The second part of the proposition is shown in Section 9.2.

9.1. *Proof of the first part of Proposition 3.* We only prove the result under (24), the proof under (25) being almost the same. By combining (43) and (44) we obtain that if $F$ is such that $RF$ is $r$th times differentiable, then for all $J \subset \{1, \ldots, n\}$ there exists a sequence $\{c_{\mathbf{i}},\ \mathbf{i} \in \mathcal{I}_r \cap J^{r+1}\}$ verifying both $c_{\mathbf{i}} \in ]\min_{j \in J} x_j, \max_{j \in J} x_j[$ and

$$-\langle \mathbf{f}, \mathbf{t}_J^* \rangle = N_J^{-1} \sum_{\mathbf{i} \in \mathcal{I}_r \cap J^{r+1}} \frac{\Lambda(F)(c_{\mathbf{i}})}{r!} \phi_{\mathbf{i}}^2(\mathbf{x}^r), \tag{70}$$

where $N_J = \mathrm{Gram}(\mathbf{1}_J, \mathbf{x}_J, \dots, \mathbf{x}_J^{r-1})\gamma_J$. Let $i^* \in J$ such that

$$\inf_{i\in J} \Lambda(F)(x_i) = \Lambda(F)(x_{i^*}).$$

We have for all $c \in ]x_J^-, x_J^+[$,

$$\Lambda(F)(c) \le \Lambda(F)(x_{i^*}) + \omega(h_J).$$

Besides, by taking $\mathbf{f} = (x_1^r/R(x_1), \dots, x_n^r/R(x_n))'$ in (44) we get that

$$\frac{1}{N_J} \sum_{\mathbf{i}\in\mathcal{I}_r\cap J^{r+1}} \phi_{\mathbf{i}}^2(\mathbf{x}^r) = \frac{\|\mathbf{x}_J^r - \Pi_{\mathcal{X}_J}\mathbf{x}_J^r\|^2}{\gamma_J}.$$

Now, by using (70) and (24) we deduce that

$$\begin{aligned}\langle \mathbf{f}, \mathbf{t}_J^*\rangle &\ge -\frac{\Lambda(F)(x_{i^*}) + \omega(h_J)}{r!}\left(\frac{1}{N_J}\sum_{\mathbf{i}\in\mathcal{I}_r\cap J^{r+1}} \phi_{\mathbf{i}}^2(\mathbf{x}^r)\right)\\ &= -(\Lambda(F)(x_{i^*}) + \omega(h_J))\frac{\|\mathbf{x}_J^r - \Pi_{\mathcal{X}_J}\mathbf{x}_J^r\|^2}{\gamma_J r!}\\ &\ge v_{\mathbf{t}_J^*}(\mathbf{f}, \beta),\end{aligned}$$

and we conclude thanks to Theorem 1.

9.2. *Proof of the second part of Proposition 3.* In order to prove this second part, we apply the first part of Proposition 3.

*Evaluation of $v_{\mathbf{t}_J^*}(\mathbf{f},\beta)$.* Let us prove that for all $J \in \bigcup_{\ell=1}^{\ell_n} \mathcal{J}^{(\ell)}$,

$$v_{\mathbf{t}_J^*}(\mathbf{f},\beta) \le \kappa^* \sqrt{\log(n)}\sigma,$$

where $\kappa^*$ depends on $\alpha, \beta, s$ and $r$ only. We use Steps 1–5 in the proof of Proposition 1. For Steps 1, 2 and 5 the proof is similar to the proof of Proposition 1.

STEP 3. *For all $\mathbf{f} = (F(x_1), \dots, F(x_n))'$ with $F^{(r)} \in \mathcal{H}_s(L)$,*

$$\frac{\|\mathbf{f} - \Pi_{V_n}\mathbf{f}\|^2}{n} \le C(s,r)L^2 n^{-2(s+r)}. \tag{71}$$

PROOF. We recall that $V_n$ is the linear space generated by

$$\{\mathbf{1}_J, \mathbf{x}_J, \dots, \mathbf{x}_J^r, J \in \mathcal{J}^{\ell_n}\}.$$

Note that the vector

$$\tilde{\mathbf{f}} = \sum_{k=1}^{\ell_n}\left(F(\bar{x}_{J_k})\mathbf{1}_{J_k} + \sum_{l=1}^{r} \frac{F^{(l)}(\bar{x}_{J_k})}{l!}(\mathbf{x}_{J_k} - \bar{x}_{J_k}\mathbf{1}_{J_k})^l\right)$$

belongs to $V_n$. Hence, using that $F^{(r)} \in \mathcal{H}_s(L)$,

$$
\begin{aligned}
&\|\mathbf{f} - \Pi_{V_n}\mathbf{f}\|^2 \\
&\quad \le \|\mathbf{f} - \tilde{\mathbf{f}}\|^2 \\
&\quad = \sum_{k=1}^{\ell_n} \sum_{i \in J_k} \left( \int_{u_1 = \bar{x}_{J_k}}^{x_i} \int_{u_2 = \bar{x}_{J_k}}^{u_1} \cdots \int_{u_r = \bar{x}_{J_k}}^{u_{r-1}} (F^{(r)}(u_r) - F^{(r)}(\bar{x}_{J_k}))\, du_r \cdots du_1 \right)^2 \\
&\quad \le \sum_{k=1}^{\ell_n} \sum_{i \in J_k} L^2 \ell_n^{-2(r+s)} \\
&\quad \le C(s,r) L^2 n^{1-2(r+s)}
\end{aligned}
$$

since $\ell_n \ge n/(4(r+1))$ using that $[x] \ge x/2$ for $x \ge 1$. $\square$

STEP 4. *Assuming that $n \ge (L/\sigma)^{1/(r+s)}$, there exists some constant $C$ depending on $s, r$ and $\beta$ only such that*

$$
\frac{\bar{\chi}^{-1}_{n-d_n, \|\mathbf{f} - \Pi_{V_n}\mathbf{f}\|^2/\sigma^2}(\beta/2)}{n - d_n} \le C. \tag{72}
$$

The proof is similar to the proof of Step 4 in Proposition 1 by using (71).

*Evaluation of $\gamma_J$.* Let us prove that there exists some constant $C$ depending on $r$ only such that, for $J$ such that $|J| \ge r+1$,

$$
\gamma_J^2 \ge C \frac{|J|^{2r+1}}{n^{2r}}.
$$

Since for all $i$, $x_i = i/n$, by translation

$$
\begin{aligned}
\gamma_J^2 &= \|\mathbf{x}_J^r - \Pi_{\mathcal{X}_J}\mathbf{x}_J^r\|^2 \\
&= \frac{1}{n^{2r}} \min_{a_0, \ldots, a_{r-1}} \sum_{i=1}^{|J|} (i^r - a_0 - a_1 i - \cdots - a_{r-1} i^{r-1})^2.
\end{aligned}
$$

By setting for all $j \in \{0, \ldots, r-1\}$ $a_j = b_j |J|^{r-j}$, we have

$$
\begin{aligned}
&\min_{a_0, \ldots, a_{r-1}} \sum_{i=1}^{|J|} (i^r - a_0 - a_1 i - \cdots - a_{r-1} i^{r-1})^2 \\
&\qquad = |J|^{2r+1} \min_{b_0, \ldots, b_{r-1}} \frac{1}{|J|} \sum_{i=1}^{|J|} \left( \left( \frac{i}{|J|} \right)^r - b_0 - \cdots - b_{r-1} \left( \frac{i}{|J|} \right)^{r-1} \right)^2.
\end{aligned}
$$

Since

$$\min_{b_0,\dots,b_{r-1}} \frac{1}{|J|}\sum_{i=1}^{|J|}\left(\left(\frac{i}{|J|}\right)^r - b_0 - \cdots - b_{r-1}\left(\frac{i}{|J|}\right)^{r-1}\right)^2$$

converges as $|J| \to \infty$ toward

$$\min_{b_0,\dots,b_{r-1}} \int_0^1 (x^r - b_0 - \cdots - b_{r-1}x^{r-1})^2\,dx,$$

which is positive, we obtain that there exists some constant $C > 0$ such that for $|J|$ large enough,

$$\gamma_J^2 \ge C\frac{|J|^{2r+1}}{n^{2r}}.$$

Moreover, since for $|J| \ge r+1$, $\gamma_J^2 > 0$, the above inequality holds for $|J| \ge r+1$, possibly enlarging $C$.

*Evaluation of* $\omega(h_J)$. Let $J \in \mathcal{J}^{(\ell)}$. Since $F^{(r)} \in \mathcal{H}_s(L)$, and since $h_J$ defined in Theorem 3 satisfies $0 < h_J \le 1/\ell$,

$$\begin{aligned}\omega(h_J) &= \sup_{|x-y|\le h_J} |F^{(r)}(x) - F^{(r)}(y)| \\ &\le L\ell^{-s}.\end{aligned}$$

*Conclusion.* Let us prove in conclusion that if

$$\inf_{x\in[0,1]} F^{(r)}(x) \le -\rho_{n,r}, \tag{73}$$

then (24) holds for some $J \in \bigcup_{\ell=1}^{\ell_n} \mathcal{J}^{(\ell)}$.

Since $F^{(r)} \in \mathcal{H}_s(L)$ under (73), there exists $j \in \{1,\dots,n\}$ such that

$$F^{(r)}(x_j) \le -\rho_{n,r} + Ln^{-s} \le -\rho_{n,r}/2$$

for $n$ large enough.

Let

$$\ell(n) = \left[\left(\frac{L^2 n}{\sigma^2\log(n)}\right)^{1/(1+2r+2s)}\right].$$

For $n$ large enough, $\ell(n) \in \{1,\dots,\ell_n\}$. Let $J$ be the element of $\mathcal{J}^{(\ell(n))}$ containing $j$. Note that $|J| \ge n/(2\ell(n))$ at least for $n$ large enough. This implies that, for $n$ large enough,

$$\gamma_J^2 \ge C\frac{|J|^{1+2r}}{n^{2r}} \ge C(r)n(\ell(n))^{-1-2r}.$$

It follows that

$$\begin{aligned} v_{\mathbf{t}_J^*}(\mathbf{f},\beta)\frac{r!}{\gamma_J}+\omega(h_J) &\le \frac{\kappa^* r!}{\sqrt{C(r)}}\sigma\sqrt{\log(n)}\frac{(\ell(n))^{r+1/2}}{\sqrt{n}}+L(\ell(n))^{-s} \\ &\le \kappa L^{(1+2r)/(1+2r+2s)}\left(\frac{\sigma^2\log(n)}{n}\right)^{s/(1+2s+2r)} \end{aligned}$$

for some constant $\kappa$ depending on $\alpha,\beta,s$ and $r$. This concludes the proof of the proposition.

## APPENDIX

### A.1. Proof of Lemma 2.

PROOF OF (42). Clearly, one has $\tilde{\mathcal{C}}_{r,R}\subset\mathcal{C}_{r,R}$. We prove $\mathcal{C}_{r,R}\subset\tilde{\mathcal{C}}_{r,R}$ by using repeatedly the following claim.

CLAIM 1. *Let $0\le u_1<u_2<\cdots<u_{r+1}<u_{r+2}\le 1$ be an increasing sequence of $r+2$ points of $[0,1]$. Let $v_1,\dots,v_{r+2}$ be real numbers verifying that*

$$D_1(1,u_{r+2},\dots,u_{r+2}^{r-1},v_{r+2})=\det\begin{pmatrix}1 & u_2 & \cdots & u_2^{r-1} & v_2\\ 1 & u_3 & \cdots & u_3^{r-1} & v_3\\ \vdots & \vdots & \vdots & \vdots & \vdots\\ 1 & u_{r+2} & \cdots & u_{r+2}^{r-1} & v_{r+2}\end{pmatrix}\ge 0$$

*and*

$$D_{r+2}=\det\begin{pmatrix}1 & u_1 & \cdots & u_1^{r-1} & v_1\\ 1 & u_2 & \cdots & u_2^{r-1} & v_2\\ \vdots & \vdots & \vdots & \vdots & \vdots\\ 1 & u_{r+1} & \cdots & u_{r+1}^{r-1} & v_{r+1}\end{pmatrix}\ge 0.$$

*Then for all $j\in\{2,\dots,r+1\}$*

$$D_j(1,u_{r+2},\dots,u_{r+2}^{r-1},v_{r+2})=\det\begin{pmatrix}1 & u_1 & \cdots & u_1^{r-1} & v_1\\ \vdots & \vdots & \vdots & \vdots & \vdots\\ 1 & u_{j-1} & \cdots & u_{j-1}^{r-1} & v_{j-1}\\ 1 & u_{j+1} & \cdots & u_{j+1}^{r-1} & v_{j+1}\\ \vdots & \vdots & \vdots & \vdots & \vdots\\ 1 & u_{r+2} & \cdots & u_{r+2}^{r-1} & v_{r+2}\end{pmatrix}\ge 0.$$

PROOF. For real numbers $t_1,\dots,t_r$ we denote by $\mathrm{vand}(t_1,\dots,t_r)$ the Vandermonde determinant

$$\mathrm{vand}(t_1,\dots,t_r)=\det\begin{pmatrix}1 & t_1 & \cdots & t_1^{r-1}\\ \vdots & \vdots & \vdots & \vdots\\ 1 & t_r & \cdots & t_r^{r-1}\end{pmatrix},$$

and for $j = 1, \dots, r+2$ we denote by $\mathbf{u}_j$ the vector $(1, u_j, \dots, u_j^{r-1}, v_j)'$. Let us fix $j \in \{2, \dots, r+1\}$. By expanding the determinant

$$D_j(1, u_{r+2}, \dots, u_{r+2}^{r-1}, v_{r+2})$$

by its last column, we get that if $j \in \{2, \dots, r\}$,

$$\begin{aligned} &D_j(1, u_{r+2}, \dots, u_{r+2}^{r-1}, v_{r+2}) \\ &\quad = v_{r+2} \operatorname{vand}(u_1, \dots, u_{j-1}, u_{j+1}, \dots, u_{r+1}) + D_j(1, u_{r+2}, \dots, u_{r+2}^{r-1}, 0), \end{aligned}$$

and if $j = r+1$,

$$\begin{aligned} &D_{r+1}(1, u_{r+2}, \dots, u_{r+2}^{r-1}, v_{r+2}) \\ &\quad = v_{r+2} \operatorname{vand}(u_1, \dots, u_r) + D_{r+1}(1, u_{r+2}, \dots, u_{r+2}^{r-1}, 0). \end{aligned}$$

Since the $u_i$'s are increasing, the Vandermonde determinants are positive and therefore $D_j(1, u_{r+2}, \dots, u_{r+2}^{r-1}, v_{r+2})$ is increasing with respect to $v_{r+2}$. On the other hand, since by assumption

$$\begin{aligned} &D_1(1, u_{r+2}, \dots, u_{r+2}^{r-1}, v_{r+2}) \\ &\quad = v_{r+2} \operatorname{vand}(u_2, \dots, u_{r+1}) + D_1(1, u_{r+2}, \dots, u_{r+2}^{r-1}, 0) \geq 0 \end{aligned}$$

we have that

$$v_{r+2} \geq -\frac{D_1(1, u_{r+2}, \dots, u_{r+2}^{r-1}, 0)}{\operatorname{vand}(u_2, \dots, u_{r+1})} = v^*,$$

and deduce that

$$D_j(1, u_{r+2}, \dots, u_{r+2}^{r-1}, v_{r+2}) \geq D_j(1, u_{r+2}, \dots, u_{r+2}^{r-1}, v^*).$$

It remains to show that $D_j(1, u_{r+2}, \dots, u_{r+2}^{r-1}, v_{r+2}, v^*) \geq 0$. When $v_{r+2} = v^*$, we have that $D_1(1, u_{r+2}, \dots, u_{r+2}^{r-1}, v^*) = 0$ and therefore $\mathbf{u}^* = (1, u_{r+2}, \dots, u_{r+2}^{r-1}, v^*)'$ is a linear combination of $\mathbf{u}_2, \dots, \mathbf{u}_{r+1}$. Let us denote by $\lambda_k$ the coordinate of $\mathbf{u}^*$ on $\mathbf{u}_k$. By Cramér's formula we have that for $k \in \{3, \dots, r\}$

$$\begin{aligned} \lambda_k &= \frac{\operatorname{vand}(u_2, \dots, u_{k-1}, u_{k+1}, \dots, u_{r+2})}{\operatorname{vand}(u_2, \dots, u_{k-1}, u_{k+1}, \dots, u_{r+1}, u_k)} \\ &= (-1)^{r-k+1} \frac{\operatorname{vand}(u_2, \dots, u_{k-1}, u_{k+1}, \dots, u_{r+2})}{\operatorname{vand}(u_2, \dots, u_{r+1})}, \end{aligned}$$

$$\lambda_2 = (-1)^{r-1} \frac{\operatorname{vand}(u_3, \dots, u_{r+2})}{\operatorname{vand}(u_2, \dots, u_{r+1})}$$

and

$$\lambda_{r+1} = \frac{\operatorname{vand}(u_2, \dots, u_r, u_{r+2})}{\operatorname{vand}(u_2, \dots, u_{r+1})}.$$

Hence, the positivity of the Vandermonde determinants implies that $\lambda_j$ has the sign of $(-1)^{r-j+1}$. Since $\mathbf{u}^* = \sum_{k=2}^{r+1} \lambda_k \mathbf{u}_k$, by linearity of the determinant

$$\begin{aligned} D_j(1, u_{r+2}, \ldots, u_{r+2}^{r-1}, v^*) &= \lambda_j D_j(1, u_j, \ldots, u_j^{r-1}, v_j) \\ &= (-1)^{r-j+1} \lambda_j D_{r+1} \end{aligned}$$

and thus, as $D_{r+1} \geq 0$, $D_j(1, u_{r+2}, \ldots, u_{r+2}^{r-1}, v^*) \geq 0$. $\square$

The proof of (42) is complete. $\square$

Proof of (43). For $x \in [x_{i_1}, x_{i_{r+1}}]$ let us set

$$h(x) = \det \begin{pmatrix} 1 & x & \cdots & x^{r-1} & R(x)F(x) \\ 1 & x_{i_2} & \cdots & x_{i_2}^{r-1} & R(x_{i_2})F(x_{i_2}) \\ \vdots & \vdots & \vdots & \vdots & \vdots \\ 1 & x_{i_{r+1}} & \cdots & x_{i_{r+1}}^{r-1} & R(x_{i_{r+1}})F(x_{i_{r+1}}) \end{pmatrix}$$

$$- \lambda \det \begin{pmatrix} 1 & x & \cdots & x^{r-1} & x^r \\ 1 & x_{i_2} & \cdots & x_{i_2}^{r-1} & x_{i_2}^r \\ \vdots & \vdots & \vdots & \vdots & \vdots \\ 1 & x_{i_{r+1}} & \cdots & x_{i_{r+1}}^{r-1} & x_{i_{r+1}}^r \end{pmatrix},$$

where $\lambda$ is such that $h(x_1) = 0$. Since $h$ is $r$-times differentiable and satisfies $h(x_{i_1}) = h(x_{i_2}) = \cdots = h(x_{ir+1}) = 0$, there exists some $c_{\mathbf{i}} \in ]x_{i_1}, x_{i_{r+1}}[$ such that

$$0 = h^{(r)}(c_{\mathbf{i}}) = \det \begin{pmatrix} 0 & 0 & \cdots & 0 & \Lambda(F)(c_{\mathbf{i}}) \\ 1 & x_{i_2} & \cdots & x_{i_2}^{r-1} & R(x_{i_2})F(x_{i_2}) \\ \vdots & \vdots & \vdots & \vdots & \vdots \\ 1 & x_{i_{r+1}} & \cdots & x_{i_{r+1}}^{r-1} & R(x_{i_{r+1}})F(x_{i_{r+1}}) \end{pmatrix}$$

$$- \lambda \det \begin{pmatrix} 0 & 0 & \cdots & 0 & r! \\ 1 & x_{i_2} & \cdots & x_{i_2}^{r-1} & x_{i_2}^r \\ \vdots & \vdots & \vdots & \vdots & \vdots \\ 1 & x_{i_{r+1}} & \cdots & x_{i_{r+1}}^{r-1} & x_{i_{r+1}}^r \end{pmatrix},$$

leading to $\lambda = \Lambda(F)(c_{\mathbf{i}})/r!$. We get the result by substituting the expression of $\lambda$ in the equality $h(x_1) = 0$. $\square$

Proof of (44). We start with the following claim.

CLAIM 2. *Let $\mathcal{W}$ be a linear subspace of $\mathbb{R}^k$ of dimension $q \in \{1, \dots, k-1\}$ and let $\{\mathbf{w}^1, \dots, \mathbf{w}^q\}$ be a basis of $\mathcal{W}$. Then for all $\mathbf{u}, \mathbf{v}$ in $\mathbb{R}^k$*

$$\begin{aligned} &\mathrm{Gram}(\mathbf{w}^1, \dots, \mathbf{w}^q)\langle \mathbf{u}, (I - \Pi_{\mathcal{W}})\mathbf{v}\rangle \\ &\quad = \sum_{\mathbf{i} \in \mathcal{I}_{q+1}} \det \begin{pmatrix} w^1_{i_1} & \cdots & w^q_{i_1} & u_{i_1} \\ \vdots & \vdots & \vdots & \vdots \\ w^1_{i_{q+1}} & \cdots & w^q_{i_{q+1}} & u_{i_{q+1}} \end{pmatrix} \\ &\qquad \times \det \begin{pmatrix} w^1_{i_1} & \cdots & w^q_{i_1} & v_{i_1} \\ \vdots & \vdots & \vdots & \vdots \\ w^1_{i_{q+1}} & \cdots & w^q_{i_{q+1}} & v_{i_{q+1}} \end{pmatrix}, \end{aligned} \tag{74}$$

*where*

$$\mathrm{Gram}(\mathbf{w}^1, \dots, \mathbf{w}^q) = \det(G) \qquad \text{with } G = (\langle \mathbf{w}^i, \mathbf{w}^j \rangle)_{1 \le i,j \le q}.$$

We conclude thanks to Claim 2 by taking $\mathbf{u} = R \star \mathbf{f}$, $\mathbf{v} = \mathbf{x}^r_J - \Pi_{\mathcal{X}_J} \mathbf{x}^r_J$, $\mathcal{W} = \mathcal{X}_J$ and $k = |J|$.

PROOF OF CLAIM 2. For $\mathbf{z} \in \mathbb{R}^k$, let $B(\mathbf{z})$ the $k \times (q+1)$ matrix

$$B(\mathbf{z}) = \begin{pmatrix} w^1_1 & \cdots & w^q_1 & z_1 \\ \vdots & \vdots & \vdots & \vdots \\ w^1_k & \cdots & w^q_k & z_k \end{pmatrix}.$$

We obtain the result by computing

$$\det(B(\mathbf{u})' B(\mathbf{v})) = \det \begin{pmatrix} \langle \mathbf{w}^1, \mathbf{w}^1 \rangle & \cdots & \langle \mathbf{w}^1, \mathbf{w}^q \rangle & \langle \mathbf{w}_1, \mathbf{v} \rangle \\ \vdots & \vdots & \vdots & \vdots \\ \langle \mathbf{w}^q, \mathbf{w}^1 \rangle & \cdots & \langle \mathbf{w}^q, \mathbf{w}^q \rangle & \langle \mathbf{w}_q, \mathbf{v} \rangle \\ \langle \mathbf{u}, \mathbf{w}^1 \rangle & \cdots & \langle \mathbf{u}, \mathbf{w}^q \rangle & \langle \mathbf{u}, \mathbf{v} \rangle \end{pmatrix}$$

by two different ways. The first way is direct: since $\Pi_W \mathbf{v}$ is a linear combination of the $\mathbf{w}^j$'s we have

$$\begin{aligned} &\det(B(\mathbf{u})' B(\mathbf{v})) \\ &\quad = \det \begin{pmatrix} \langle \mathbf{w}^1, \mathbf{w}^1 \rangle & \cdots & \langle \mathbf{w}^1, \mathbf{w}^q \rangle & \langle \mathbf{w}^1, (I - \Pi_W)\mathbf{v} \rangle \\ \vdots & \vdots & \vdots & \vdots \\ \langle \mathbf{w}^q, \mathbf{w}^1 \rangle & \cdots & \langle \mathbf{w}^q, \mathbf{w}^q \rangle & \langle \mathbf{w}^q, (I - \Pi_W)\mathbf{v} \rangle \\ \langle \mathbf{u}, \mathbf{w}^1 \rangle & \cdots & \langle \mathbf{u}, \mathbf{w}^q \rangle & \langle \mathbf{u}, (I - \Pi_W)\mathbf{v} \rangle \end{pmatrix} \\ &\quad = \det \begin{pmatrix} \langle \mathbf{w}^1, \mathbf{w}^1 \rangle & \cdots & \langle \mathbf{w}^1, \mathbf{w}^q \rangle & 0 \\ \vdots & \vdots & \vdots & \vdots \\ \langle \mathbf{w}^q, \mathbf{w}^1 \rangle & \cdots & \langle \mathbf{w}^q, \mathbf{w}^q \rangle & 0 \\ \langle \mathbf{u}, \mathbf{w}^1 \rangle & \cdots & \langle \mathbf{u}, \mathbf{w}^q \rangle & \langle \mathbf{u}, (I - \Pi_W)\mathbf{v} \rangle \end{pmatrix} \end{aligned}$$

$$= \mathrm{Gram}(\mathbf{w}^1, \ldots, \mathbf{w}^q)\langle \mathbf{u}, (I - \Pi_W)\mathbf{v}\rangle.$$

The other way is to use the Cauchy–Binet formula [see Horn and Johnson (1991)]: we calculate $\det(B(\mathbf{u})'B(\mathbf{v}))$ as a function of the $(q+1)\times(q+1)$ minors of the matrix $B(\mathbf{u})$ and $B(\mathbf{v})$ which leads to the right-hand side of (74) and concludes the proof. □

The proof of (44) is complete. □

### A.2. Proof of Proposition 2.

*Case* $\mathcal{K} = \mathcal{K}_{\geq 0}$. Let $\mathbb{P}_F$ be the law of $Y$ under the model defined by (16). Let $\Phi$ be a test of level $\alpha$ of the hypothesis $F \in \mathcal{K}_{\geq 0}$. Let us define the test $\Psi$ of the hypothesis "$F = 0$" against "$F \neq 0$" which rejects the null if $\Phi(Y) = 1$ or if $\Phi(-Y) = 1$. Since $0 \in \mathcal{K}_{\geq 0}$ and since

$$\mathbb{P}_0(\Phi(Y) = 1) = \mathbb{P}_0(\Phi(-Y) = 1) \leq \alpha,$$

the test $\Psi$ is of level $2\alpha \leq 3\alpha$. Let $\rho_n(\Phi, \mathcal{F})$ be the $\Delta$-uniform separation rate of $\Phi$ over $\mathcal{F}$. It is enough to show that

$$\rho_n(\Phi, \mathcal{F}) \geq \rho_n(0, \mathcal{F}).$$

To do so, we show that the $\|\cdot\|_\infty$-uniform separation rate of $\Psi$ over $\mathcal{F}$ is not larger than $\rho_n(\Phi, \mathcal{F})$, which means that for all $F \in \mathcal{F}$ such that $\|F\|_\infty \geq \rho_n(\Phi, \mathcal{F})$ we have $\mathbb{P}_F(\Psi(Y) = 1) \geq 1 - \beta$.

Let $F \in \mathcal{F}$. If $\|F\|_\infty \geq \rho_n(\Phi, \mathcal{F})$, then

$$\text{either} \quad \Delta(F) = \sup_{x \in [0,1]} (-F(x)\mathbf{1}_{F(x)>0}) \geq \rho_n(\Phi, \mathcal{F}) \quad \text{or} \quad \Delta(-F) \geq \rho_n(\Phi, \mathcal{F}).$$

In the first case, by definition of $\rho_n(\Phi, \mathcal{F})$ we have $\mathbb{P}_F(\Phi(Y) = 1) \geq 1 - \beta$ and consequently $\mathbb{P}_F(\Psi(Y) = 1) \geq 1 - \beta$. Note that in the other case the same is true since by symmetry of the law of $Y - F$

$$\mathbb{P}_F(\Phi(-Y) = 1) = \mathbb{P}_{-F}(\Phi(Y) = 1).$$

*Case* $\mathcal{K} = \mathcal{K}_{\nearrow}$. We argue similarly. Let $\Phi$ be a test of level $\alpha$ of the hypothesis $F \in \mathcal{K}_{\nearrow}$. We also consider the test $\Phi'$ of level $\alpha$ of "$F = 0$" against "$F \neq 0$" which rejects the null when $\sqrt{n}|\int_0^1 dY(t)|$ is large enough (namely, larger than the $1 - \alpha$ quantile of a standard Gaussian random variable). Finally, we define the test $\Psi$ of the hypothesis "$F = 0$" against "$F \neq 0$" which rejects the null if $\Phi(Y) = 1$ or $\Phi(-Y) = 1$ or $\Phi'(Y) = 1$. Since $0 \in \mathcal{K}_{\nearrow}$, we have that the so-defined test $\Psi$ is of level $3\alpha$.

Some easy computations show that there exists some constant $\kappa$ depending on $\alpha$ and $\beta$ only such that $\Phi'$ rejects the null with probability not smaller than $1 - \beta$ as soon as $|\int_0^1 F(t)\, dt|$ is larger than $\kappa\sigma/\sqrt{n}$ (the sum of the $\beta$ and $1-\alpha$ quantiles of a standard Gaussian suits for $\kappa$). On the other hand, note that

$$\Delta(F) = \frac{1}{2} \sup_{0 \leq s \leq t \leq 1} (F(s) - F(t))$$

and thus, by definition of the $\Delta$-separation rate, $\rho_n(\Phi, \mathcal{F})$, of $\Phi$ over $\mathcal{F}$, $\Psi$ rejects the null with probability not smaller than $1-\beta$ under all alternatives $F \in \mathcal{F}$ satisfying

$$\max\{\Delta(F), \Delta(-F)\} = \frac{1}{2} \sup_{0 \leq t, s \leq 1} |F(t) - F(s)| > \rho_n(\Phi, \mathcal{F}).$$

Therefore, since

$$\begin{aligned} \|F\|_\infty &\leq \sup_{t \in [0,1]} \left| F(t) - \int_0^1 F(s)\,ds \right| + \left| \int_0^1 F(s)\,ds \right| \\ &\leq \int_0^1 \sup_{t \in [0,1]} |F(t) - F(s)|\,ds + \left| \int_0^1 F(s)\,ds \right| \\ &\leq \sup_{t,s \in [0,1]} |F(t) - F(s)| + \left| \int_0^1 F(s)\,ds \right|, \end{aligned}$$

$\Psi$ rejects the null with probability larger than $1-\beta$ under all alternative $\mathcal{F}$ such that

$$\|F\|_\infty \geq 2\rho_n(\Phi, \mathcal{F}) + \kappa\sigma/\sqrt{n},$$

and the result follows.

**Acknowledgment.** The authors would like to thank Claudine Picaronny for pointing out to us the proof given for (74).

Y. BARAUD
LABORATOIRE J. A. DIEUDONNE
UMR CNRS 6621
UNIVERSITÉ DE NICE SOPHIA-ANTIPOLIS
PARC VALROSE
06108 NICE CEDEX 2
FRANCE
E-MAIL: baraud@math.unice.fr

S. HUET
INRA, LABORATOIRE DE BIOMÉTRIE
78352 JOUY-EN-JOSAS CEDEX
FRANCE
E-MAIL: sylvie.huet@jouy.inra.fr

B. LAURENT
INSA, DEPARTEMENT GMM
135 AV. DE RANGUEIL
31077 TOULOUSE CEDEX 4
FRANCE
E-MAIL: beatrice.laurent@insa-toulouse.fr